\documentclass[10pt,twocolumn,twoside,journal]{IEEEtran} 

   \usepackage[pdftex]{graphicx}
   \DeclareGraphicsExtensions{.pdf}

\usepackage{cite}
\usepackage[hyperfootnotes=false]{hyperref}
\hypersetup{ 		
    pdfauthor={Florian Dorfler},        
    colorlinks=true,        
    linkcolor=black,
    urlcolor=black,
    citecolor=black     		
}

\usepackage[usenames,dvipsnames]{color}
\usepackage{epstopdf}
\usepackage{subfigure}

\usepackage{amsmath,amssymb,mathrsfs,dsfont,mathdots}

\usepackage[nocomma]{optidef}

\newtheorem{theorem}{Theorem}[section]
\newtheorem{lemma}[theorem]{Lemma}
\newtheorem{fact}[theorem]{Fact}
\newtheorem{definition}[theorem]{Definition}
\newtheorem{corollary}[theorem]{Corollary}
\newtheorem{proposition}[theorem]{Proposition}

\newtheorem{remark}[theorem]{Remark}

\newcommand\oprocendsymbol{\hbox{$\square$}}
\newcommand\oprocend{\relax\ifmmode\else\unskip\hfill\fi\oprocendsymbol}

\newcommand{\real}[0]{\mathbb R}

\DeclareSymbolFont{bbold}{U}{bbold}{m}{n}
\DeclareSymbolFontAlphabet{\mathbbold}{bbold}

\definecolor{jgreen}{rgb}{0.0, 0.5, 0.0}

\newcommand{\tb}[0]{\color{black}}

\usepackage[prependcaption,colorinlistoftodos]{todonotes}

\newcommand{\bv}{{\mathscr{B}}}
\newcommand{\Han}{{\mathscr{H}}}
\newcommand{\lti}{{\mathscr{L}}}
\newcommand{\obs}{{\mathscr{O}}}
\newcommand{\conv}{{\mathscr{G}}}
\newcommand{\ctr}{{c_\text{ctrl}}}
\newcommand{\cid}{{c_\text{id}}}
\newcommand{\Tini}{{T_{\textup{ini}}}}

\newcommand{\Td}{T}
\newcommand{\Tr}{L}
\newcommand{\uini}{{u_{\textup{ini}}}}
\newcommand{\yini}{{y_{\textup{ini}}}}
\newcommand{\wini}{{w_{\textup{ini}}}}
\newcommand{\alphani}{{x_{\textup{ini}}}}
\newcommand{\Up}{{U_{\mathrm{p}}}}
\newcommand{\Uf}{{U_{\mathrm{f}}}}
\newcommand{\Yp}{{Y_{\mathrm{p}}}}
\newcommand{\Yf}{{Y_{\mathrm{f}}}}

\usepackage{amsmath}

\DeclareMathOperator*{\argmin}{arg\,min}
\DeclareMathOperator*{\where}{where}
\DeclareMathOperator*{\st}{subject\,\,to}

\usepackage{accents}
\newlength{\dhatheight}

\title{
Bridging direct \& indirect data-driven control formulations via regularizations and relaxations
}

\author{Florian D\"orfler, Jeremy Coulson, and Ivan Markovsky
\thanks{FD and JC are with Department of Information Technology and Electrical Engineering, ETH Zurich, 8092 Zurich, Switzerland. Email: {\tt\small \{dorfler,coulson\}@control.ee.ethz.ch}. IM is with the Department of Electrical Engineering (ELEC) of the Vrije Universiteit Brussel, 1050 Brussels, Belgium. Email: {\tt\small imarkovs@vub.be}.}
\thanks{This work was supported by ETH Zurich funds. IM was supported by the Fond for Scientific Research Vlaanderen (FWO) project G090117N and the Excellence of Science Project 30468160.}}

\begin{document}

\maketitle
\thispagestyle{empty}
\pagestyle{empty}

\begin{abstract}
We discuss connections between sequential system identification and control for linear time-invariant systems, 
{\tb often termed}
indirect data-driven control, as well as a  {\tb contemporary} direct data-driven control approach seeking an optimal decision compatible with recorded data assembled in a Hankel matrix and robustified through suitable regularizations. We formulate these two problems in the language of behavioral systems theory and parametric mathematical programs, and we bridge them through a multi-criteria formulation trading off system identification and control objectives. We illustrate our results with two methods from subspace identification and control: namely, subspace predictive control and low-rank approximation which constrain trajectories to be consistent with a non-parametric predictor derived from (respectively, the column span of) a data Hankel matrix. In both cases we conclude that direct and regularized data-driven control can be derived as convex relaxation of the indirect approach, and the regularizations account for an implicit identification step. Our analysis further reveals a novel regularizer and {\tb a plausible   hypothesis explaining the remarkable empirical} performance of direct methods on nonlinear systems.%
\end{abstract}



\section{Introduction}
\label{Sec: Intro}

The vast realm of data-driven control methods can be classified into {\em indirect data-driven control} approaches consisting of sequential system identification and model-based control as well as {\em direct data-driven control} approaches seeking an optimal decision compatible with recorded data. Both approaches have a rich history, and they have received renewed interest cross-fertilized by novel methods and widespread interest in machine learning. Representative recent surveys are \cite{hewing2020learning,pillonetto2014kernel,chiuso2019system,hou2013model,recht2019tour,IM-FD:21-survey,hjalmarsson2005experiment,SM:22}.

The pros and cons of both paradigms have often been elaborated on. {\tb The indirect approach is modular with well understood subtasks, though modeling and identification are cumbersome, their results are often not useful for control (due to, e.g., incompatible uncertainty quantifications), and practitioners  often prefer end-to-end methods.  Direct approaches promise to resolve these problems by learning control policies directly from data. However, they are often analytically and computationally less tractable and rarely apply to real-time and safety-critical control systems. Selected direct methods that proved themselves in theory and practice are iterative feedback tuning and virtual reference feedback tuning~\cite{HH-MG-SG-OL:98, MC-AL-SS:02,bazanella2011data}.

 Quite a few approaches have bridged the direct and indirect data-driven control paradigms.} Of relevance to this article, we note the literature on identification for control \cite{hjalmarsson2005experiment,hjalmarsson1996,geversaa2005,schrama1992} and control-oriented regularized identification \cite{formentin2018core}, which propose that the control objective should bias the identification task. {\tb Likewise, dual control dating to \cite{feldbaum1963dual} addresses the  {exploration vs. exploitation} trade-offs in simultaneous identification and optimal control; see \cite{ferizbegovic2019learning,larsson2016application,iannelli2020structured} for  recent contributions. Furthermore, \cite{campestrini2017data} formulates data-driven model reference control as an identification problem, where various degrees of prior information can be incorporated so that the method can range between the direct and the indirect approach.}

We take a similar perspective here:  the sequential identification and control tasks can be abstracted as nested bi-level optimization problem: find the best control subject to a model, where the model is the best fit to a data set within some hypothesis class. This approach is modular and both steps admit tractable formulations, but generally it is also suboptimal: there is no separation principle -- aside from special cases, see \cite[Section 4]{hjalmarsson2005experiment} -- for these two nested optimization problems. An end-to-end direct algorithmic approach may thus outperform indirect methods if a tractable formulation were available. 
For the latter we resort to a paradigm square in between behavioral system theory and subspace system identification methods.

Behavioral system theory \cite{willems2007,willems1991,willems1997} takes an abstract view on dynamical systems as sets of trajectories, and it does not require parametric representations which makes it appealing from a data-centric perspective. For example, linear time-invariant (LTI) systems are characterized as shift-invariant subspaces within an ambient  space of time series. The role of identification is to find such a low-dimensional feature from data. Subspace methods take a similar (albeit more algorithmic) viewpoint \cite{van2012subspace,katayama2006subspace,van1994n4sid} and extract parametric models from the range and null spaces of a low-rank data Hankel matrix.

Both lines of work come together in a result known as the Fundamental Lemma \cite{willems2005}; see also \cite{waarde2020,IM-FD:20,IM-FD:21-survey} for recent extensions. It states that, under some assumptions, the set of all finite-length trajectories (the restricted behavior) of an LTI system equals the range space of a data Hankel matrix. This result serves as the theoretic underpinning for work in subspace identification \cite{IM-FD:20,markovsky2006,markovsky2005algorithms}
and data-driven control, in particular subspace predictive control based on non-parametric models \cite{favoreel1999spc,qin2005novel,huang2008dynamic},  explicit feedback policies parametrized by data matrices \cite{berberich2020combining,van2020noisy,de2019formulas},
and data-enabled predictive control {\tb (DeePC)} seeking compatibility of predicted trajectories with the range space of a data Hankel matrix. The latter methods have first been established for deterministic LTI systems in \cite{markovsky2008,markovsky2016} and have recently been extended by suitably regularizing the optimal control problems. Closed-loop stability was certified in \cite{berberich2020data}. The regularizations were first mere heuristics \cite{JC-JL-FD:18} but have later been constructively derived by robust control and optimization \cite{JC-JL-FD:19-CDC,JC-JL-FD:20,xue2020data,LH-JZ-JL-FD:20,LH-JZ-JL-FD:01}.  These approaches, albeit recent, have proved themselves in practical nonlinear problems {\tb in multiple domains} \cite{LH-JZ-JL-FD:01,LH-JC-JL-FD:19,PC-AF-SB-FD:20,EE-JC-PB-JL-FD:19,LH-JZ-JL-FD:20}. We also note the recent maximum-likelihood perspective \cite{yin2020maximum}.
{\tb We refer to \cite{IM-FD:21-survey} surveying results surrounding the fundamental lemma}.

In this paper, we explore the following questions: how does data-enabled predictive control relate to a prior system identification? What are principled regularizations? And why does it work so well in the nonlinear case? We start our investigations from indirect data-driven control formulated as a bi-level optimization problem  in the general output feedback setting. 
As a vehicle to transition between indirect and direct approaches, we consider a multi-criteria problem trading off identification and control objectives {\tb reminiscent of similar approaches \cite{hjalmarsson2005experiment,hjalmarsson1996,geversaa2005,schrama1992,formentin2018core,feldbaum1963dual,ferizbegovic2019learning,larsson2016application,iannelli2020structured,campestrini2017data} blending the two.} We formally show that one tail of its Pareto front corresponds to the bi-level problem, and a convex relaxation results in the regularized data-enabled predictive control formulations used in \cite{berberich2020data,JC-JL-FD:18,JC-JL-FD:19-CDC,JC-JL-FD:20,xue2020data,LH-JZ-JL-FD:20,LH-JC-JL-FD:19,PC-AF-SB-FD:20,EE-JC-PB-JL-FD:19,LH-JZ-JL-FD:01}. 

Most of our results are formulated in the abstract language of behavioral systems theory and parametric mathematical programs, but we also specialize our treatment to two concrete methods: subspace predictive control {\tb (SPC)} \cite{favoreel1999spc,qin2005novel,huang2008dynamic} and low-rank approximation \cite{markovsky2016}. In both cases we conclude that the direct regularized data-driven control can be derived as a convex relaxation of the indirect approach, where {\tb $(i)$ LTI complexity specifications (selecting the model class) are dropped, and $(ii)$} the projection of the data on the set of LTI systems is replaced by regularizations accounting for implicit identification. 
{\tb In particular, starting from indirect data-driven control based on low-rank approximation of a Hankel matrix, we arrive at a DeePC formulation with an $\ell_{1}$-regularizer (Theorem \ref{Theorem: Low-rank relaxation}). When formulating indirect data-driven control via the SPC framework, our analysis reveals a novel regularizer for DeePC promoting a least-square data fit  by projecting on the null space of the Hankel matrix (Theorem~\ref{Theorem: SPC relaxation}).}

We illustrate our results with numerical studies illustrating the role of regularization, superiority of the new regularizer, and comparisons. Informed by our analysis, we hypothesize and numerically confirm that the indirect approach is superior in case of ``variance'' error, e.g., for LTI stochastic systems, and the direct approach wins in terms of ``bias'' error, e.g., for nonlinear systems supporting the {\tb empirical} observations in \cite{LH-JZ-JL-FD:01,LH-JC-JL-FD:19,PC-AF-SB-FD:20,EE-JC-PB-JL-FD:19,LH-JZ-JL-FD:20}. 
{\tb Similar bias-variance trade-offs can also be found in the recent pre-print \cite{krishnan2021direct} discussing sub-optimality of direct and indirect methods as  function of the data size.\,These findings also resonate with those of data-driven model reference control \cite{campestrini2017data} concluding that the direct approach is superior in reducing the bias whereas the indirect one gives better variance -- especially if an erroneous model class is selected.}

The remainder of this paper is organized as follows: Section~\ref{Sec: Preliminaries} reviews representations of LTI systems. Section~\ref{Sec: Direct and Indirect Data-Driven Control} formulates the direct and indirect data-driven control problems, and Section~\ref{Sec: Bridging} bridges them. Section~\ref{subsec: numerical analysis} contains our numerical studies. Finally, Section~\ref{sec: conclusions} concludes the paper. {\tb Readers familiar with the behavioral approach may skip Section \ref{Sec: Preliminaries}.}


\section{LTI Systems and their Representations}
\label{Sec: Preliminaries}

We adopt a behavioral perspective which allows for system theory independent of parametric representations.
We aim at a concise exposition and refer to \cite{willems2007,willems1991,willems1997,IM-FD:21-survey} for details.

\subsection{Behavioral Perspective on Discrete-Time LTI systems}

Consider the discrete time axis $\mathbb Z$, the signal space $\real^{q}$, and the associated {space of trajectories}  $\real^{q\mathbb Z}$ consisting of all $q$-variate sequences $(\dots,w(-1),w(0), w(1),\dots)$ with $w(i) \in \real^{q}$. 
Consider a permutation matrix $P$ partitioning each 
$
w(i)
= P
\left[\begin{smallmatrix}
u(i) \\ y(i)
\end{smallmatrix}\right]
$,
where $u(i) \in \real^{m}$ and  $y(i) \in \real^{q-m}$ are free and dependent variables that will later serve as inputs and outputs. The {\em behavior} $\bv$ is defined as a subset of the space of trajectories, $\bv \subset \real^{q\mathbb Z}$, and a system as the triple $(\mathbb Z,\real^{q},\bv)$. 

In what follows, we denote a system merely by its behavior $\bv$, keeping the  signal space  $ \real^{q\mathbb Z}$ fixed throughout. A system is {\em linear} if $\bv$ is a subspace of $ \real^{q\mathbb Z}$. Let $\sigma$ denote the shift operator with action $\sigma w({t}) = w({t+1})$. A system is {\em time-invariant} if $\bv$ is shift-invariant: $\sigma \bv = \bv$. 
%
{\tb Finally, $\bv_{L}$ is the {restriction of $\bv$} to $ \real^{qL}$, i.e., to trajectories of length $L\in \mathbb Z_{>0}$.}


\subsection{Kernel Representations and Parametric Models}
\label{subsec: kernel reps}

{\tb 
Rather than a mere set-theoretic descriptions, one typically works with explicit {\em parametric representations} (colloquially termed {\em models}) of LTI systems. 
For instance, a {\em kernel representation} with {\em lag} $\ell$ specifies an LTI behavior  as}
\begin{equation*}
\bv = \text{kernel}(R(\sigma)) = \bigl\{ w \in \real^{q\mathbb Z}\,:\; R(\sigma) w = 0 \bigr\}\,,
\end{equation*}
where $R(\sigma) = R_{0}+R_{1} \sigma + \dots + R_{\ell} \sigma^{\ell}$ is a polynomial matrix of degree $\ell$, and the matrices $R_{0}, R_{1}, \dots, R_{\ell}$ take values in $\real^{(q-m) \times q}$. 
Alternatively, one can unfold the kernel representation  by revealing a latent variable: the state $x(t) \in \real^{n}$. The {\em input/state/output} (or {\em state-space}) {\em representation} is 
\begin{align*}
	\bv = & \bigl\{  w = P \left[\begin{smallmatrix} u \\ y\end{smallmatrix}\right]  \in \real^{q\mathbb Z}\,:\; \exists x \in \real^{n\mathbb Z} \mbox{ such that }
	\\
	&  \quad \sigma x = Ax + Bu\,,\, y = Cx+Du \bigr\}\,,
\end{align*}
where $A \in \real^{n \times n}$,  $B \in \real^{n \times m}$, $C \in \real^{q-m \times n}$, and  $D \in \real^{q-m \times m}$.
We assume that the lag $\ell$ (resp., the state dimension $n$) is minimal, i.e., there is no other kernel (resp., state-space) representation with smaller lag (resp., state dimension).   

The dimension $n$ of a minimal state-space representation manifests itself in a minimal kernel representation as $n=\sum_{i=1}^{q-m} \ell_{i}$, where $\ell_{i}$ is the {degree}\ 
of the $i$th row of $R(\sigma)$.

\subsection{Representation-Free Estimation and Behavior Dimension}
\label{subsec: estimation}

{\tb
Given a state-space representation with $m$ inputs, order $n$, and lag $\ell$, the extended {observability} and convolution matrices
\begin{equation*}
\obs_{L}
=
\left[\begin{smallmatrix}
C \\ CA \\ \vdots \\ C A^{L-1}
\end{smallmatrix}\right]
\quad\mbox{and}\quad
\conv_{L}
=
\left[\begin{smallmatrix}
D & 0 & \cdots & & 0
\\
CB & D & 0 & \cdots & 0
\\
CAB & CB & D & \ddots & \vdots
\\
\vdots & \ddots & \ddots &\ddots & 0
\\
CA^{L-2}B & \cdots & CAB & CB & D
\end{smallmatrix}\right]
\end{equation*}
 parametrize all length-$L$ trajectories in $\bv_{L}$ as
\begin{equation}
\begin{bmatrix}
u \\ y
\end{bmatrix}
=
\begin{bmatrix}
I & 0 \\ \conv_{L} & \obs_{L}
\end{bmatrix}
\begin{bmatrix}
 u \\ \alphani 
\end{bmatrix}
\,,
\label{eq: IOS representation}
\end{equation}
where $\alphani \in \real^{n} $ is the initial state.
Recall the {\em observability problem}: given length-$L$ time series of inputs and outputs, can $\alphani$ be reconstructed?  Equation \eqref{eq: IOS representation} gives a succinct answer: namely, $\alphani$ can be reconstructed if and only if $\obs_{L}$ has full column-rank. The minimum $L$ so that $\obs_{L}$ has full rank $n$ equals the {\em lag} $\ell$ of a minimal kernel representation. 
As readily deducible from \eqref{eq: IOS representation} and formalized in \cite[Lemma\,1]{markovsky2008}, in a {\em representation-free} setting, the  initial condition $\alphani$ for a trajectory $w \in \bv_{L}$, can be estimated via a  prefix trajectory $\wini = \bigl(w(-\Tini+1),  \dots, w(-1), w(0) \bigr)$  of length $\Tini \geq \ell$ so that the concatenation $\wini \wedge w \in \bv_{\Tini+L}$ is a valid trajectory.}

Hence, an LTI system is characterized by the complexity parameters $(	q,m,n,\ell)$, and we denote the corresponding  class of LTI systems by $\lti_{m,\ell}^{q,n}$: namely, LTI systems with $m$ inputs, $q-m$ outputs, minimal state dimension $n$, and minimal lag $\ell$.

The following lemma characterizes  the dimension of $\bv_{L} \in \lti_{m,\ell}^{q,n}$ in terms of the complexity parameters $(	q,m,n,\ell)$.

{\tb
\begin{lemma}[Dimension of $\bv_{L}$]
\label{Lemma: subspace dimension}
Let $\bv \in \lti_{m,\ell}^{q,n}$. Then $\bv_{L}$ is a subspace of $ \real^{qL}$, and for $L \geq \ell$ its  dimension is $mL + n$.
\end{lemma}

\begin{IEEEproof}
Due to linearity of $\bv$, $\bv_{L} \subset  \real^{qL}$ is a subspace. To show that the dimension of $\bv_{L}$ equals $mL + n$ for $L \geq \ell$, we appeal to a minimal state-space representation of $\bv$ --- a state-space-independent proof is in \cite[Sec. 3]{IM-FD:20}. We have  $w = P \left[\begin{smallmatrix} u \\ y\end{smallmatrix}\right] \in \bv_{L}$ if and only if \eqref{eq: IOS representation}  holds
for some $\alphani \in \real^{n}$. Since the representation is minimal,  $\obs_{L} \in \real^{(q-m)L \times n}$ is of full column-rank for $L \geq \ell$. Therefore, the matrix $
\left[\begin{smallmatrix}
I & 0 \\ \conv_{L} & \obs_{L}
\end{smallmatrix}\right]
\in \real^{qL \times (mL + n)}$
is of full rank $mL+n$ for  $L \geq \ell$ and forms a basis for $\bv_{L}$. Thus, $\bv_{L}$ has dimension $mL+n$. 
\end{IEEEproof}

\begin{remark}[Complexity bounds]
All forthcoming results assume known complexity $(q,m,n,\ell)$.
When only data and no prior information is available, it is reasonable to assume upper bounds on $(q,m,n,\ell)$. In this case, the anticipated dimension of $\bv_{L}$ is at most $mL + n$, and the forthcoming rank equalities the behavior dimension should be replaced by  inequalities.
\oprocend
\end{remark}}

\subsection{Image Representation of Restricted Behavior}\label{sec:image-representation}

The restricted behavior $\bv_{L}$, the set of all trajectories of length $L$, can be described by a kernel or state-space representation. As an interesting alternative, we recall  the  {\em image representation} of $\bv_{L}$ by a data matrix of a time series.

Consider the sequence $w = \bigl(w(1),w(2),\dots,w(T)\bigr)$ with elements  $w(i) \in \real^{q}$, and define the (block) {\em Hankel matrix} $\Han_{L}(w) 
\in \real^{qL \times (T-L+1)}$ of depth $L$, for some $L \leq T$, as
\begin{equation*}
\Han_{L}(w) =
\begin{bmatrix}
w(1) & w(2) & \dots & w(T-L+1)
\\
w(2)  & w(3) & \dots & w(T-L+2)
\\
\vdots & \vdots  & \ddots & \vdots
\\
w(L) & w(L+1) & \dots & w(T)
\end{bmatrix}
\,.
\end{equation*}
A result due to \cite{willems2005} that became known as the {\em Fundamental Lemma}  offers an image representation of the restricted behavior in terms of the column span of a data Hankel matrix. We present a necessary and sufficient version here assuming:
\begin{enumerate}\addtolength{\itemindent}{5pt}
\renewcommand{\labelenumi}{{\theenumi}}\renewcommand{\theenumi}{(A.\arabic{enumi})}
\item \label{ass:L+n pe} rank$\left(\Han_L(w)\right)=mL + n$. 
\end{enumerate}

\begin{lemma}\label{lemma: fundamental lemma}{\cite[Corollary 19]{IM-FD:20})}:
Consider an LTI system $\bv \in \lti_{m,\ell}^{q,n}$ and an associated trajectory  $w = \bigl(w(1), w(2),$ $ \dots, w(T)\bigr) \in \real^{qT}$. The following are equivalent for $L > \ell$:\vspace{-4pt}
\begin{equation*}
	\text{colspan} \left(\Han_{L}(w) \right) = \bv_{L}  
	\quad\Longleftrightarrow\quad
	\text{Assumption \ref{ass:L+n pe}}
\end{equation*}
\end{lemma}

In words, the Hankel matrix $\Han_L(w)$ composed of a single $T$-length  trajectory parametrizes all $L$-length trajectories if and only if rank$\left(\Han_L(w)\right)=mL + n$. A plausible reasoning leading up to Lemma~\ref{lemma: fundamental lemma} is that every column of $\Han_L(w)$ is a trajectory of length $L$, and the set of all such trajectories has at most dimension $mL+n$; see Lemma~\ref{Lemma: subspace dimension}.
Lemma~\ref{lemma: fundamental lemma} extends the original {\em Fundamental Lemma}  \cite[Theorem 1]{willems2005} which requires input/output partitioning, controllability, and persistency of excitation of order $L+n$ (i.e., $\Han_{L+n}(u)$ must have full row rank) as sufficient conditions. Lemma~\ref{lemma: fundamental lemma} also extends to mosaic Hankel, Page, and trajectory matrices \cite{{IM-FD:20}}.

{\tb
\begin{remark}[Models vs. data]\label{rem: models vs. data}
 It is debatable whether the image representation via the Hankel matrix $\Han_{L}(w)$ should be called a ``model'', as it is readily available from raw data. Hence, we call $\text{colspan} \left(\Han_{L}(w) \right)$ a {\em data-driven representation} of $\bv_{L}$ and reserve the term ``model'' for parametric (kernel or state-space) representations. Models are useful  for many reasons: first and foremost the availability of powerful analysis and design methods. Another readily discernible advantage is that models are vastly compressed compared to the image representation, and the latter holds only on finite horizons unless trajectories are weaved together \cite{markovsky2005algorithms}; see also Remark~\ref{rem: On data lengths}.
\oprocend
\end{remark}}


\section{Direct and Indirect Data-Driven Control}
\label{Sec: Direct and Indirect Data-Driven Control}

We present different data-driven control formulations along with assumptions under which the formulations are consistent. These assumptions are used only for consistency statements and not for our main results, but they will prove insightful.


\subsection{Optimal Control Problem}
\label{subsec:opt-ctr}

Given a plant with {\em plant behavior} $\bv^P \in \lti_{m,\ell}^{q,n}$, a $\Tini$-length prefix trajectory $\wini = \bigl(w(-\Tini+1),  \dots, w(0) \bigr) \in \bv_{\Tini}$, a $\Tr$-length reference trajectory $w_r\in\real^{q\Tr}$ in a {\em reference behavior} $\bv^R$,
{\tb and a set of {\em admissible trajectories\/} $\mathcal W \subset \real^{q\Tr}$,}\ 
consider the finite-time {\em optimal control problem}
\begin{mini}
{w {\tb \in \mathcal W}}{\ctr(w-w_{r}) }{\label{eq:OPT}}{\boldsymbol{C}:}
\addConstraint{
{\wini \wedge w \in \bv^{P}_{\Tini+\Tr}}
\,.
}
\end{mini}
For  $\Tini \geq \ell$ the prefix trajectory $\wini$ implicitly sets the initial condition for the optimal control problem \eqref{eq:OPT}; see Section~\ref{subsec: estimation}. In case of uncertain initial condition, the prefix $\wini$ can be made a decision variable and included via a penalty term in the cost; c.f., \cite{berberich2020data,markovsky2016,JC-JL-FD:18,JC-JL-FD:20,JC-JL-FD:19-CDC}. We refrain from such extensions here.%

{\tb Typically, the cost $\ctr:\, \mathbb R^{q\Tr} \to \real_{\geq 0}$ includes a running and a terminal cost. The set $\mathcal W \subset \real^{q\Tr}$ captures constraints  on admissible trajectories (e.g., capturing input saturation). We denote a minimizer (if it exists) of the optimization problem $\boldsymbol{C}$ in \eqref{eq:OPT} by $w^{\star}_{C}$.
We make the following regularity assumptions:
\begin{enumerate}\addtolength{\itemindent}{5pt}
\renewcommand{\labelenumi}{{\theenumi}}\renewcommand{\theenumi}{(A.\arabic{enumi})}\addtocounter{enumi}{1}


\item\label{ass:ctr cost} $\ctr:\, \mathbb R^{q\Tr} \to \real_{\geq 0}$ is a convex function that achieves its minimum when $w = w_{r}$;  $\mathcal W \subset \real^{q\Tr}$ is closed,\,convex, and non-empty; and $(\real^{q\Tini} \oplus \mathcal W) \cap \bv^{P}_{\Tini+\Tr}$ is\,{non-empty.}

  
	
\end{enumerate}
The last assumption ensures that $\mathcal W$ is {\em viable}, i.e., a trajectory of $\bv^{P}$ originating anywhere can be contained within $\mathcal W$ for $L$ steps.
Problem \eqref{eq:OPT} is thus convex with closed, convex, and non-empty feasible set due to Assumption~\ref{ass:ctr cost} and because ${\bv^{P}_{\Tini+\Tr}}$ is a subspace; see Lemma \ref{Lemma: subspace dimension}. Under further standard assumptions existence and uniqueness of a (global) minimum can be assured, but we do not impose further structure.
 
 For problem \eqref{eq:OPT}, we do not necessarily assume $\bv^P = \bv^R$, since we often ask systems to track non-plant behavior (e.g., steps). Likewise, we  generally do not assume feasibility: $w_{r} \in \mathcal W$.  However, such assumptions connect to {model reference control} and allow to state  consistency results as presented next. 
\begin{enumerate}\addtolength{\itemindent}{5pt}
\renewcommand{\labelenumi}{{\theenumi}}\renewcommand{\theenumi}{(A.\arabic{enumi})}\addtocounter{enumi}{2}

\item\label{ass: Bp = Br} $\wini \wedge w_{r} \in {(\real^{q\Tini} \oplus \mathcal W)} \cap \bv^{P}_{\Tini+\Tr}$, i.e., the reference $w_{r} \in \bv^R_{\Tr}$ is compatible with  the prefix trajectory $\wini$, the plant $\bv^{P}$, and the constraints $\mathcal W$.
	
\end{enumerate}
}

\begin{fact}\label{Fact: minimum of C}
Under Assumptions \ref{ass:ctr cost} and \ref{ass: Bp = Br}, the minimum of the  control problem $\boldsymbol{C}$ in \eqref{eq:OPT} is   achieved for $w^{\star}_{C}=w_{r}$.

\end{fact}

{\tb Fact~\ref{Fact: minimum of C} (and similar consistency results later) follows since $w^{\star}_{C}=w_{r}$ is feasible and achieves the minimum of the cost.
Fact~\ref{Fact: minimum of C} (and consistency Assumption \ref{ass: Bp = Br}) serve to establish ground-truth for comparing different problem formulations.}

Problem \eqref{eq:OPT} becomes a ``classical'' control problem if  a parametric model for the plant $\bv^P$ is available. The latter is usually obtained from data through system identification. 

\subsection{Indirect Data-Driven Control via System Identification}
\label{subsec:id+ctr}

Given a $\Td$-length trajectory $w_{d} \in\real^{q\Td}$ as  {\em identification data},  conventional system identification and control consists of three steps. The first step, {\em model class selection}, amounts to choosing the set of candidate models, e.g., $\lti_{m,\ell}^{q,n}$ specified by the complexity $(q,n,m,\ell)$.
The second step, {\em model fitting}, chooses an element from the model class that fits the data best in some specified sense, e.g., distance between data $w_{d}$ and model $\bv$.
This step is often synonymous to learning a parametric model (e.g., PEM), though some classic (e.g., ETFE) and modern (e.g., kernel-based) methods are non-parametric and by-pass the model {\tb order} selection; see \cite{pillonetto2014kernel} for a review (and the acronyms). However, for control design the non-parametric models again have to be projected on a behavior in $\lti_{m,\ell}^{q,n}$. 
Both approaches can be abstracted\,as
\begin{mini}
{\hat w_{d},\widehat\bv}{\cid(\hat w_{d}-w_{d}) }{\label{eq:ID}}{\boldsymbol{ID}:}
\addConstraint{
\hat w_{d} \in \widehat\bv_{\Td}\,,\; \widehat\bv \in \lti_{m,\ell}^{q,n}\,.
}
\end{mini}
It is useful to think of the identification loss $\cid: \mathbb R^{q\Td} \!\to\! \real_{\geq 0}$ as a distance. Given the data $w_{d}$, problem \eqref{eq:ID} seeks the closest 
LTI behavior  within the class $\lti_{m,\ell}^{q,n}$, i.e., the closest subspace with dimensions as in Lemma \ref{Lemma: subspace dimension}. 
We denote a minimizer of \eqref{eq:ID} by $\bigl(\hat w^{\star}_{d,ID},\widehat\bv_{ID}^{\star}\bigr)$ and assume the following {about}\ 
$\cid(\cdot)$:
\begin{enumerate}\addtolength{\itemindent}{5pt}
\renewcommand{\labelenumi}{{\theenumi}}\renewcommand{\theenumi}{(A.\arabic{enumi})}\addtocounter{enumi}{3}

\item \label{ass: cid} $\cid(\cdot)$ achieves its minimum when $\hat w_{d} \!=\! w_{d}$.
	
	\end{enumerate}
{\tb Note that existence and uniqueness of minimizers of \eqref{eq:ID} does not only hinge upon the regularity of cost and constraint functions, but also on the data. In general, identification problems are non-convex. For now we keep problem \eqref{eq:ID} abstract and general and resort to more specific formulations in Section~\ref{Sec: Bridging}.}

Exact identification of the true system requires exact data $w_{d} \in \bv^P_{\Td}$ and an identifiability assumption \cite[Theorem 15]{IM-FD:20} which assures that $\bv^P$ can be  recovered from $w_{d}$:
\begin{enumerate}\addtolength{\itemindent}{5pt}
\renewcommand{\labelenumi}{{\theenumi}}\renewcommand{\theenumi}{(A.\arabic{enumi})}\addtocounter{enumi}{4}

	\item \label{ass: Bp=Bid} $w_{d}   \in \bv^P_{\Td}$, i.e., $w_{d}$ is a valid trajectory of $\bv^P_{\Td}$	; and

	\item \label{ass: ell + n + 1 pe} rank$\left(\Han_{\ell+1}(w_{d})\right)=m(\ell+1) + n$.

\end{enumerate}

\begin{fact}\label{Fact: minimum of ID}
Under assumptions \ref{ass: cid}--\ref{ass: ell + n + 1 pe}, the minimum value of the system identification problem $\boldsymbol{ID}$ in \eqref{eq:ID} is achieved for $\hat w_{d,ID}^{\star}=w_{d}$ and $\widehat\bv^{\star}_{ID} = \bv^P$.
\end{fact}

{\tb We again note that the (arguably strong) Assumptions  \ref{ass: Bp = Br}, \ref{ass: Bp=Bid}, and \ref{ass: ell + n + 1 pe} are used only for consistency statements (such as Fact~\ref{Fact: minimum of ID}) and not for our later main results and simulations.}

Finally, equipped with an identified behavior $\widehat\bv^{\star} \in \lti_{m,\ell}^{q,n}$, the third step is {\em certainty-equivalence control}: solve the optimal control problem \eqref{eq:OPT} subject to the identified model:
\begin{mini}
{ w {\tb \in \mathcal W}}{\ctr( w-w_{r}) }{\label{eq:OPT-CE}}{}
\addConstraint{
{\wini \wedge  w  \in \widehat\bv^{\star}_{\Tini+\Tr}}
\,.
}
\end{mini}
In \eqref{eq:OPT-CE}, $\ctr( w-w_{r})$ is merely a {\em surrogate} (predicted) control error since  $ w \in \widehat\bv^{\star}$, the identified model, rather than $ w \in \bv^{P}$. 

Putting both the system identification \eqref{eq:ID} and  certainty-equivalence control \eqref{eq:OPT-CE} together, we arrive at  {indirect data-driven control} formulated as the {\em bi-level problem}%
\begin{mini}%
{w {\tb \in \mathcal W}}{\ctr(w-w_{r}) }{\label{eq:OPT-BL}}{\!\!\!\!\boldsymbol{BL}:\!\!\!}
\addConstraint{
{\wini \wedge w  \in \widehat\bv^{\star}_{\Tini+\Tr}}
}
\addConstraint{\where \quad
\widehat\bv^{\star} \in \argmin_{\hat w_{d},\,\,\,\widehat{\bv}} \; \cid(\hat w_{d}-w_{d})
}
\addConstraint{\qquad\qquad\! \st \quad \hat w_{d} \in\,\widehat\bv_{\Td}\,,
}
\addConstraint{\qquad\qquad \qquad\qquad\quad\;\, \widehat\bv \in \lti_{m,\ell}^{q,n}
\,.
}
\end{mini}
The bi-level problem structure in \eqref{eq:OPT-BL} reflects the sequential system identification and control tasks, that is, first a model is fitted to the data in the inner  identification problem before the model is used for control in the outer problem.
We denote a minimizer for  the inner problem of \eqref{eq:OPT-BL} by $\bigl(\hat w^{\star}_{d,BL},\widehat \bv^{\star}_{BL}\bigr)$ and a minimizer for the outer problem of \eqref{eq:OPT-BL} by $w^{\star}_{BL}$.

{\tb
\begin{remark}[Further problem levels and the value of models]
The bi-level formulation \eqref{eq:OPT-BL} is only the tip of the iceberg, and the overall design may feature further nested levels, e.g., optimization of the model selection hyper-parameters $(n,\ell)$, uncertainty quantification, etc. We deliberately neglect these levels here and focus on identification and control.

Since our ultimate interest is control, we treat models in a disregarding manner, i.e., they serve merely an auxiliary purpose. Of course, models are desired for other reasons: system design, analysis, the reasons in Remark~\ref{rem: models vs. data}, etc.
\oprocend
\end{remark}}

Under suitable consistency assumptions, the sequential system identification and control approach in \eqref{eq:OPT-BL} is optimal.

\begin{fact}\label{Fact: minimum of BL}
Consider the optimal control problem $\boldsymbol{C}$ in \eqref{eq:OPT} and the bi-level problem $\boldsymbol{BL}$ in \eqref{eq:OPT-BL}. Then
\begin{enumerate}

\item under Assumptions \ref{ass: cid}--\ref{ass: ell + n + 1 pe}, the bi-level problem $\boldsymbol{BL}$ reduces to the optimal control $\boldsymbol{C}$; and

\item under the additional Assumptions \ref{ass:ctr cost} and \ref{ass: Bp = Br}, the minimum value of the  bi-level problem $\boldsymbol{BL}$ is achieved for $\hat w_{d,BL}^{\star}=w_{d}$, and $\widehat\bv^{\star}_{BL} = \bv^P$, $w^{\star}_{BL}= w_{r}$.
\end{enumerate}
\end{fact}

The first statement echos the ``model as well as possible'' paradigm and a separation of control\,and\,identification, albeit in a simple setting; see \cite[Section 4.2]{hjalmarsson2005experiment} for further reading. 

\subsection{Direct Data-Driven Control via the Image Representation}
\label{subsec: Direct Data-Driven Control via the Image Representation}

The direct data-driven control approach pursued here hinges upon the Fundamental Lemma~\ref{lemma: fundamental lemma}. A direct corollary of {\tb the latter} is that the prediction and estimation trajectories have to be within the column span of the data Hankel matrix.

\begin{corollary}[Direct data-driven control]\label{Corollary: Data-driven optimal control}
Assume that Assumptions \ref{ass:L+n pe} and \ref{ass: Bp=Bid} hold with $L$ replaced by $\Tini+L$, then the optimal control problem $\boldsymbol{C}$ in \eqref{eq:OPT} is equivalent to
\begin{mini}
{w {\tb\in \mathcal W}}{\hspace{-7pt}\ctr(w-w_{r}) }{\label{eq:OPT-D}}{\hspace{-30pt}\boldsymbol{D}\!:\!\!\!}
\addConstraint{\hspace{-7pt}
\begin{bmatrix}
\wini\\w\end{bmatrix}
\!\in\! \text{colspan} \left(\Han_{\Tini+\Tr}(w_{d}) \right)\!,
}
\end{mini}
i.e., the minimizers and minima of \eqref{eq:OPT} and \eqref{eq:OPT-D} coincide.
\end{corollary}

\begin{fact}
Under Assumptions \ref{ass:L+n pe} with $L$ replaced by $\Tini+L$, \ref{ass:ctr cost}, \ref{ass: Bp = Br}, and \ref{ass: Bp=Bid}, the minimum value of \eqref{eq:OPT-D} is  achieved for $w^{\star}_{D} = w_{r}$.
\end{fact}

{\tb
\begin{remark}[Data lengths]\label{rem: On data lengths}
It is instructive to compare the sample complexity of direct and indirect approaches \eqref{eq:OPT-D} and \eqref{eq:OPT-BL}.
Due to Assumption \ref{ass:L+n pe},  \eqref{eq:OPT-D} requires more data than the identification Assumption {\ref{ass: ell + n + 1 pe}}. This discrepancy is due to \eqref{eq:OPT-D} seeking a multi-step predictor, whereas identification \eqref{eq:ID} seeks a single-step predictor to be applied recursively.  
  By weaving multiple trajectories of length $\ell+1$, Assumption \ref{ass:L+n pe} can be eased so that the data lengths coincide;   see \cite[Lemma 3]{markovsky2005algorithms}. 
\oprocend
\end{remark}

In comparison, with system identification, the model order selection is implicit in Assumption \ref{ass:L+n pe} and encoded in the rank of the Hankel matrix $\Han_{\Tini+\Tr}(w_{d})$ -- at least, for exact data $w_{d} \in \bv^{P}_{\Td}$. If the data $w_{d}$ is noisy, then $\Han_{\Tini+\Tr}(w_{d})$ likely has full rank, and the constraint  of \eqref{eq:OPT-D} is vacuous. Thus, $w=w_{r}$ uniquely minimizes the surrogate control error, but the realized control error may be arbitrarily different. In short, certainty-equivalence can fail arbitrarily poorly in direct data-driven control, and the direct approach has to be robustified. This is a major difference with the indirect (first identify, then control) approach \eqref{eq:OPT-BL}: one purpose of identification is to filter noisy data by projecting on a deterministic behavior.

To go beyond certainty equivalence, the DeePC approaches \cite{berberich2020data,JC-JL-FD:19-CDC,JC-JL-FD:20,xue2020data,JC-JL-FD:18,LH-JZ-JL-FD:20,LH-JZ-JL-FD:01} reformulate the constraint in \eqref{eq:OPT-D} as $\text{col}(\wini,w) = \Han_{\Tini+\Tr}(w_{d})  g$ for some $g$ and add a robustifying regularizer.}%
\begin{mini}
{w {\tb \in \mathcal W},g}{\ctr(w-w_{r}) + \lambda \cdot h(g)}{\label{eq:OPT-DR}}{\boldsymbol{D}_{\lambda}:}
\addConstraint{
\begin{bmatrix}
\wini\\w\end{bmatrix} = \Han_{\Tini+\Tr}(w_{d})  g
}
\end{mini}
To provide an intuition, every column of $\Han_{\Tini+\Tr}(w_{d})$ is a trajectory of $\bv^{P}_{\Tini+\Tr}$, and the decision variable $g$ linearly combines these columns for the optimal trajectory $w$ -- consistent with the prefix trajectory $\wini$ and  regularized by $h(g)$.
The regularization function $h(\cdot)$ and parameter $\lambda$ are nonnegative. Choices for $h(\cdot)$ are one-norms \cite{JC-JL-FD:18}, two-norms \cite{xue2020data}, squared two-norms \cite{berberich2020data,LH-JZ-JL-FD:20}, or arbitrary $p$-norms \cite{JC-JL-FD:19-CDC,JC-JL-FD:20,LH-JZ-JL-FD:01}. 

{\tb
The regularizers can be related to robust optimization formulations in  deterministic \cite{xue2020data,LH-JZ-JL-FD:20,LH-JZ-JL-FD:01} or stochastic settings \cite{JC-JL-FD:19-CDC,JC-JL-FD:20}, where $\lambda$ is a design parameter specifying  the size of the assumed uncertainty set. The regularized formulation \eqref{eq:OPT-DR} has proved itself in  practical (nonlinear) control systems \cite{LH-JZ-JL-FD:01,LH-JC-JL-FD:19,PC-AF-SB-FD:20,EE-JC-PB-JL-FD:19,LH-JZ-JL-FD:20}.}


\section{Bridging Direct \& Indirect Approaches}
\label{Sec: Bridging}

\subsection{Multi-Objective Data-Driven Control}
\label{subsec: multi-objective}

From an optimization perspective it is natural to lift the bi-level problem \eqref{eq:OPT-BL} to a {\em multi-criteria problem} simultaneously optimizing for identification and control objectives. Using weighted sum scalarization, the multi-criteria problem\,is\vspace{-5pt}%
\begin{mini}%
{w {\tb \in \mathcal W},\hat w_{d},\widehat\bv}{\gamma \cdot \cid(\hat w_{d}-w_{d}) \,+\, \ctr(w-w_{r}) }{\label{eq:OPT-SMO}}{\!\!\!\boldsymbol{MC}_{\gamma}\!:\!\!\!}
\addConstraint{
{\wini \wedge  w  \in \widehat\bv_{\Tini+\Tr}}\,,\; \hat w_{d} \in \widehat\bv_{\Td}\,,}
\addConstraint{ \widehat\bv \in  \lti_{m,\ell}^{q,n}
\,,
}
\end{mini}
where the trade-off parameter $\gamma \geq 0$ traces the Pareto front between the identification and optimal control objectives. 

The multi-criteria problem \eqref{eq:OPT-SMO} can be interpreted as fitting a model $\widehat\bv$ simultaneously to two data sets: the identification data $w_{d}$ and the  reference $w_{r}$. From a control perspective, the identification criterion biases the solution $w \in \widehat\bv$ to adhere to the observed data $w_{d}$ rather than merely matching the to be tracked reference $w_{r}$. Likewise, from the other side, the identification criterion is biased by the control objective.
In short, control and identification {\em regularize} each other, in the spirit of identification for control \cite{hjalmarsson2005experiment,hjalmarsson1996,geversaa2005,schrama1992}. A  similar  formulation has been proposed in \cite{formentin2018core} interpolating between PEM identification and a model-reference control objective. {\tb Likewise, the data-driven model reference control formulation in \cite{campestrini2017data} interpolates between a direct and an indirect approach. Finally, dual control approaches consider similar multi-criteria formulations balancing exploration (for identification) and exploitation (i.e., optimal control) \cite{feldbaum1963dual,ferizbegovic2019learning,larsson2016application,iannelli2020structured}.}

We denote a minimizer of  \eqref{eq:OPT-SMO} by $\bigl(w^{\star}_{MC},\hat w^{\star}_{d,MC},\widehat \bv^{\star}_{MC}\bigr)$.

\begin{fact}
Under Assumptions \ref{ass:ctr cost}--\ref{ass: ell + n + 1 pe}, for any $\gamma \geq 0$ the minimum of the parametric multi-criteria problem $\boldsymbol{MC}_{\gamma}$ is  achieved for $\hat w_{d,MC}^{\star}=w_{d}$, $\widehat\bv^{\star}_{MC} = \bv^P$, and $w^{\star}_{MC} = w_{r}$.
\end{fact}

Different points on the Pareto front of \eqref{eq:OPT-SMO} have different emphasis regarding the control and identification objectives. Below we formalize that for $\gamma$ sufficiently large, the multi-criteria problem \eqref{eq:OPT-SMO} recovers the bi-level problem \eqref{eq:OPT-BL} corresponding to sequential system identification and control.

We follow standard penalty arguments from bi-level optimization \cite{ye1995,ye1997exact}, which are particularly tractable here since  \eqref{eq:OPT-BL} is only weakly coupled: the inner problem does not depend on the decision variable $w$ of the outer problem. {\tb Assume there is a minimum (termed value function) of the inner problem:}
\begin{mini}
{\hat w_{d},\,\,\,\widehat\bv}{\cid(\hat w_{d}-w_{d}) }{\label{eq:ID-value}}{\varphi \, = }
\addConstraint{
\hat w_{d} \in\,\, \widehat\bv_{\Td}\,,\;\,  \widehat\bv \in \lti_{m,\ell}^{q,n}\,.
}
\end{mini}
The bi-level problem \eqref{eq:OPT-BL} reads then equivalently as%
\begin{mini}
{w{\tb \in \mathcal W},\hat w_{d},\widehat\bv}{\ctr(w-w_{r}) }{\label{eq:OPT-BL-value}}{}
\addConstraint{
{\wini \wedge w  \in \widehat\bv_{\Tini+\Tr}}\,,\; \hat w_{d} \in \widehat\bv_{\Td}\,,
}
\addConstraint{
\widehat\bv \in \lti_{m,\ell}^{q,n}\,,\; \cid(\hat w_{d}-w_{d}) - \varphi = 0 \,.
}
\end{mini}
At this point the reader is encouraged to review the definition and salient properties of a constraint qualification termed partial calmness \cite{ye1995,ye1997exact}; see the appendix. If problem~\eqref{eq:OPT-BL-value} is partially calm at a local minimizer and $\ctr(\cdot)$ is continuous, then there is $\gamma^{\star}>0$ so that, for all $\gamma > \gamma^{\star}$, then \eqref{eq:OPT-BL-value} equals
\begin{mini}
{w{\tb \in \mathcal W},\hat w_{d},\widehat\bv}{\!\!\!\gamma \cdot \bigl| \cid(\hat w_{d}-w_{d}) - \varphi \bigr| \,+\, \ctr(w-w_{r}) }{\label{eq:OPT-BL-value-2}}{\!\!\!\!\!\!}
\addConstraint{
{\!\!\!\wini \wedge w  \in \widehat\bv_{\Tini+\Tr}}\,,\; \hat w_{d} \in \widehat\bv_{\Td}\,,\;  \widehat\bv \in \lti_{m,\ell}^{q,n} \,,
}
\end{mini}
that is, the local minimizers of \eqref{eq:OPT-BL-value} and \eqref{eq:OPT-BL-value-2} coincide; see Proposition~\ref{Proposition: partial calmness and exact penalty}.
 We now drop the absolute value (since $ \cid(\hat w_{d}-w_{d}) - \varphi \geq 0$) and the constant $\varphi$ (which in our case does not depend on the variable $w$ of the outer problem) from the objective of \eqref{eq:OPT-BL-value-2} to recover problem \eqref{eq:OPT-SMO}. 
 We have thus established a chain of equivalences relating the bi-level and multi-criteria problems. 
 We summarize our discussion below.

\begin{proposition}[Upper tail of the Pareto front of $\boldsymbol{MC}_{\gamma}$]\label{Proposition: upper tail of the Pareto front}
Consider the parametric multi-criteria problem $\boldsymbol{MC}_{\gamma}$ in \eqref{eq:OPT-SMO} and  the bi-level problem $\boldsymbol{BL}$ in \eqref{eq:OPT-BL}. Assume that {\tb the inner identification problem admits a minimum as in \eqref{eq:ID-value}, \eqref{eq:OPT-BL-value} is partially calm at any local minimizer, and $\ctr(\cdot)$ is continuous.} Then there is $\gamma^{\star}>0$ so that for $\gamma > \gamma^{\star}$ the problem $\boldsymbol{MC}_{\gamma}$ is equivalent to $\boldsymbol{BL}$, i.e., $ w^{\star}_{MC} =  w^{\star}_{BL}$ , $\hat w^{\star}_{d,MC} = \hat w^{\star}_{d,BL}$, and $\widehat \bv^{\star}_{MC} = \widehat \bv^{\star}_{BL}$.
  Moreover, the optimal values of $\boldsymbol{MC}_{\gamma}$ and $\boldsymbol{BL}$ coincide up to the constant $\gamma \cdot \varphi$ with $\varphi$ defined in \eqref{eq:ID-value}.
\end{proposition}

The following comments are in order regarding partial calmness. As discussed in Proposition \ref{Proposition: partial calmness and exact penalty}, partial calmness is equivalent to the constraint $ \cid(\hat w_{d}-w_{d}) - \varphi \geq 0$ serving as an exact penalty. Partial calmness is satisfied, for instance, appealing to Proposition~\ref{Proposition: LipschitzPenalty}, if the identification cost $\cid(\cdot)$ can be phrased as a distance (see the discussion following the identification problem \eqref{eq:ID}) and $\ctr(\cdot)$ is Lipschitz continuous over the feasible set, e.g., the feasible set is either compact (due to constraints) or the control performance is measured by a norm or Huber loss. The Lipschitz constant then serves as a lower estimate for $\gamma^{\star}$. 
%
A non-Lipschitz cost requires $\gamma \to \infty$ as a sufficient condition. Note that for $\gamma \to \infty$ Proposition~\ref{Proposition: upper tail of the Pareto front} holds without assumptions, since \eqref{eq:OPT-BL-value-2} is merely an indicator function reformulation of  \eqref{eq:OPT-BL-value}. Our relaxations in the next sections will, among others, drop the requirement on $\gamma$ sufficiently large  as well as the LTI complexity specification $\widehat\bv \in \lti_{m,\ell}^{q,n}$. 

{\tb Even if the identification \eqref{eq:ID} is convex the multi-criteria problem \eqref{eq:OPT-SMO} is not,} since it simultaneously optimizes over the to-be-identified model $\bv$ and  the to-be-designed trajectory $w$. This can be spotted in a kernel representation: the constraint $ w \in \widehat\bv_{\Tr}$ takes the form $\widehat R(\sigma) w = 0$, where both $\widehat R$ and $ w$ are variables. Other representations lead to the same conclusions.

\begin{proposition}\label{Proposition: }
Consider the multi-criteria problem \eqref{eq:OPT-SMO} and a kernel representation of the to-be-identified behavior: $\widehat\bv = \text{kernel}(\widehat R(\sigma))$. Then the feasible set of \eqref{eq:OPT-SMO}  is not convex.
\end{proposition}

We believe that the multi-criteria problem is interesting in its own right: studying its Pareto front and choosing an optimal trade-off parameter may possibly yield superior  performance.  

Our problem setup thus far was conceptual rather than practically useful. Below, we consider concrete problem formulations and turn our conceptual insights into concise results.

\subsection{Bridging Towards Subspace Predictive Control (SPC)}
\label{subsec: subspace ARX formulations}

We explain SPC  from the perspective of the Fundamental Lemma~\ref{lemma: fundamental lemma} stating that any trajectory $\wini \wedge w \in \bv_{\Tini+\Tr}^{P}$ lies in $\text{colspan} \left(\Han_{\Tini+\Tr}(w_{d}) \right)$. Recall that $\wini$ is a prefix trajectory of length $\Tini \geq \ell$ setting the initial condition, and $w$ is a future trajectory  of length $\Tr>1$ to be designed via optimal control. Accordingly, permute and partition $w$ and the Hankel matrix
\begin{equation*}
\begin{bmatrix} \wini \\ w \end{bmatrix} \sim \begin{bmatrix}
\uini \\ u \\ \hline \yini  \\ y
\end{bmatrix}
,\,\;
\Han_{\Tini+L}(w_{d})
\sim
\begin{bmatrix}
\Up \\ \Uf \\\hline \Yp  \\ \Yf
\end{bmatrix}
=
\begin{bmatrix}
\Han_{\Tini + \Tr}(u_{d}) \\\hline \Han_{\Tini + \Tr}(y_{d})
\end{bmatrix}
\,,
\end{equation*}
where $\uini \in \real^{m\Tini}$, $\yini \in \real^{(q-m)\Tini}$, and $\sim$ denotes similarity under a coordinate  permutation.
The subscripts ``p'' and ``f'' are synonymous to ``past'' and ``future''. 
We seek a linear model, i.e., a matrix $K$, relating past and future as 
\begin{equation}
	y = 
	\underbrace{
	\begin{bmatrix}
	K_{p} & \vline& K_{f}
	\end{bmatrix}}_{=K}
	\cdot \begin{bmatrix}
\uini \\ \yini \\\hline u 
\end{bmatrix}
\label{eq: ARX transition model}
\,.
\end{equation}
The multi-step predictor $K$ is found from  Hankel matrix data by means of the least-square criterion  \cite[Section\,3.4]{huang2008dynamic}
\begin{mini}
{K}{	\left\|\Yf  -  K \cdot \begin{bmatrix}
\Up \\ \Yp \\ \Uf 
\end{bmatrix}\right\|_{F}^{2}}
{\label{eq:pem}}{}
\,,
\end{mini}
where $\|\cdot\|_{F}$ is the Frobenius norm. 
Via the Moore-Penrose inverse, the solution of \eqref{eq:pem} is the classic SPC predictor \cite{favoreel1999spc}%
\begin{equation}
K =
\Yf 
\begin{bmatrix}
\Up \\ \Yp \\\Uf
\end{bmatrix}^{\dagger}
\,.
\label{eq: SPC predictor}
\end{equation}
It is insightful to compare equation \eqref{eq: ARX transition model} and the matrices $K_{p},K_{f}$ to equation \eqref{eq: IOS representation} and the extended observability and impulse response matrices $\obs_{L}$ and $\conv_{L}$, respectively. One realizes that for exact data, \eqref{eq: ARX transition model} is an ARX model with rank$(K_{p})=n$ assuring LTI behavior of desired complexity and a lower block-triangular zero pattern of $K_{f}$ assuring causality. {\tb For inexact data,  LTI behavior of desired complexity is promoted by low-rank approximation (typically via singular-value thresholding of $K_{p}$) \cite{favoreel1999spc}; and one aims to gain causality by heuristically thresholding $K_{f}$ towards a desired zero pattern
\cite[Remark 10.1]{huang2008dynamic}, \cite[Section 3]{qin2005novel}. The causality requirement can also be omitted for offline or receding horizon control, but it is useful to condition the data on the set of causal models. These steps bring the linear relation \eqref{eq: ARX transition model} half-way towards an LTI model. Though a model has further structure, e.g., $K_{f}$ is Toeplitz, and the entries of $K_{p}$ and $K_{f}$ are coupled; see \eqref{eq: IOS representation}.}

Hence, in this case the identification problem \eqref{eq:ID} {\tb is relaxed to the single, monolithic, and non-convex program}
  \begin{mini*}
{K}{	\left\|\Yf  -  K \cdot \begin{bmatrix}
\Up \\ \Yp \\ \Uf 
\end{bmatrix}\right\|_{F}^{2}}
{\label{eq:pem}}{}
\addConstraint{K = \begin{bmatrix}	 K_{p} & \vline&  K_{f} \end{bmatrix}}
\addConstraint{K_{f}\, \text{lower-block triangular}}
\addConstraint{\text{rank}( K_{p}) = n}%
\,,%
\end{mini*}%
where the lower-block triangular specification means that all entries above the diagonal $(q-m)\times m$ blocks equal zero.

 {\tb We obtain a parametric version of the indirect data-driven approach \eqref{eq:OPT-BL}, where $\wini \wedge w  \in \widehat\bv^{\star}_{\Tini+\Tr}$ and $w \in \mathcal W = \mathcal U \times \mathcal Y$ are replaced by \eqref{eq: ARX transition model} and $(u,y) \in \mathcal U \times \mathcal Y$, respectively:}%
\begin{mini}
{\tb u \in \mathcal U,y \in \mathcal Y}{\ctr\left( \left[\begin{smallmatrix} y - y_{r} \\ u - u_{r}   \end{smallmatrix} \right] \right) }{\!\!\!\!}{\label{eq:OPT-BL-ARX}}
\addConstraint{
	y = K^{\star} \cdot \begin{bmatrix}
\uini \\ \yini \\ u 
\end{bmatrix}
}
\addConstraint{\where \;
K^{\star} \in \argmin_{K} \; 
\left\|\Yf  -   K \cdot \begin{bmatrix}
\Up \\ \Yp \\ \Uf 
\end{bmatrix}\right\|_{F}^{2}
}
\addConstraint{\qquad\quad\!\; \st \quad\! \! K = \begin{bmatrix}	 K_{p} & \vline&  K_{f} \end{bmatrix}}
\addConstraint{\qquad\quad\!\;  \phantom{\st} \quad\!  K_{f}\, \text{lower-block triangular}}
\addConstraint{\qquad\quad\!\;  \phantom{\st} \quad\! \text{rank}( K_{p}) = n}
\,.
\end{mini}%
{\tb We stress that \eqref{eq:OPT-BL-ARX} is generally not an equivalent reformulation of \eqref{eq:OPT-BL} since the inner identification does not necessarily lead to an LTI model; see the comments following equation \eqref{eq: SPC predictor}.}

For comparison, consider  also an instance of the direct regularized problem \eqref{eq:OPT-DR} with regularizer $h(g) = \|(I-\Pi)g\|_{p}$:
\begin{mini}%
{\tb u \in \mathcal U,y \in \mathcal Y,g}{\ctr\left( \left[\begin{smallmatrix} y - y_{r} \\ u - u_{r}   \end{smallmatrix} \right] \right) \,+\, \lambda \cdot \|(I-\Pi)g\|_{p}  }{\label{eq:ARX-SMO-final}}{}
\addConstraint{	
\begin{bmatrix}
\Up \\ \Yp \\\Uf \\ \Yf
\end{bmatrix}
 g = \begin{bmatrix}
\uini \\ \yini \\ u \\y
\end{bmatrix} 
}\,.
\end{mini}%
Here, {\tb $\|\cdot\|_{p}$ is any $p$-norm,} $\Pi = 
\left[\begin{smallmatrix}
\Up \\ \Yp \\\Uf 
\end{smallmatrix}\right]^{\dagger}\!
\left[\begin{smallmatrix}
\Up \\ \Yp \\\Uf 
\end{smallmatrix}\right]$, and $(I-\Pi)$ is an orthogonal\,projector on the kernel of the first three block-constraint equations.
{\tb The proof of Theorem~\ref{Theorem: SPC relaxation} will later show that this regularizer is in fact  {\em induced} by the least-square identification~\eqref{eq:pem}, i.e., $\|(I-\Pi)g\|_{p}=0$ if and only if the least-square criterion is minimized. Hence, it robustifies the problem akin to least squares.} We state the following consistency result.%
\begin{fact}\label{fact: consistency of projection}
Under Assumptions \ref{ass:L+n pe} with $L$ replaced by $\Tini+L$, \ref{ass:ctr cost}, \ref{ass: Bp = Br}, and \ref{ass: Bp=Bid}, for any $\lambda \geq 0$ the minimum of the regularized  problem \eqref{eq:ARX-SMO-final} is achieved for $y^{\star} = Y_{f}g^{\star} = y_{r}$ and $u^{\star} = U_{f}g^{\star} = u_{r}$, where $\|(I-\Pi)g^{\star}\|_{p}=0$.
\end{fact}

{\tb \begin{remark}[Consistency of regularizers]\label{rem: consistency}
Fact~\ref{fact: consistency of projection} may not appear insightful at first glance, but it highlights an important fact. The projection-based regularizer $h(g)= \|(I-\Pi)g\|_{p}$ is consistent since it penalizes only the homogenous solution to the constraint equations \eqref{eq:ARX-SMO-final} and does not affect the variables $(u,y)$. In comparison,
the conventional norm-based regularizer $h(g) = \|g\|_{p}$ is {not} consistent: it penalizes the heterogeneous solution of the constraint equations in \eqref{eq:ARX-SMO-final} and thus also $(u,y)$. Hence, even with ideal consistency Assumptions \ref{ass:L+n pe}, \ref{ass:ctr cost}, \ref{ass: Bp = Br}, and \ref{ass: Bp=Bid} in place, the norm-based regularizer $h(g) = \|g\|_{p}$ with $\lambda \neq 0$ does not lead to the ground-truth solution $y^{\star} = Y_{f}g^{\star} = y_{r}$, $u^{\star} = U_{f}g^{\star} = u_{r}$; see also Remark~\ref{rem: comments on spc theorem}. \oprocend
\end{remark}}

The following is the main result of this subsection.

\begin{theorem}[SPC relaxation]\label{Theorem: SPC relaxation}
Consider the indirect data-driven control problem \eqref{eq:OPT-BL-ARX} and the direct data-driven control problem \eqref{eq:ARX-SMO-final} parameterized by $\lambda \geq 0$. {\tb Let Assumption~\ref{ass:ctr cost} hold and} assume that $\ctr(\cdot)$ is Lipschitz continuous. For $\lambda$ sufficiently small,  \eqref{eq:ARX-SMO-final} is a convex relaxation of  \eqref{eq:OPT-BL-ARX}, that is, 
\begin{enumerate}
\item[$(i)$] \eqref{eq:ARX-SMO-final} is convex,
\item[$(ii)$]  any feasible $(u,y)$ in \eqref{eq:OPT-BL-ARX}  is feasible for \eqref{eq:ARX-SMO-final}, and 
\item[$(iii)$]  the optimal value of \eqref{eq:ARX-SMO-final} lower-bounds that of \eqref{eq:OPT-BL-ARX}.
\end{enumerate}
\end{theorem}

\begin{IEEEproof}
First, we perform a convex relaxation by dropping the rank and block-triangularity constraints in \eqref{eq:OPT-BL-ARX}. Second, observe that the explicit solution of the inner problem, the predictor \eqref{eq: SPC predictor}, is equivalently derived as least-norm\,solution
\begin{align*}
	y = \Yf g^{\star}
	\;\text{ where}\;
	&g^{\star} = \argmin_{g} \| g\|_{2}
	\nonumber
	\\&\text{subject to}
		\begin{bmatrix}
\Up \\ \Yp \\\Uf 
\end{bmatrix}
 g = \begin{bmatrix}
\uini \\ \yini \\ u 
\end{bmatrix} 
\,.
\end{align*}
\mbox{We now insert this reformulation in the relaxation of \eqref{eq:OPT-BL-ARX}:}\vspace{-0pt}%
\begin{mini}%
{\tb u \in \mathcal U,y \in \mathcal Y}{\ctr\left( \left[\begin{smallmatrix} y - y_{r} \\ u - u_{r}   \end{smallmatrix} \right] \right) }{\!\!\!\!\!\!\!\!\!\!\!\!\!\!\!\!\!\!\!\!\!\!\!\!\!\!\!}{\label{eq:least-norm form}}
\addConstraint{	y = \Yf g^{\star}}
\addConstraint{\where \quad
g^{\star} \in \argmin_{ g} \; \| g\|_{2}}
\addConstraint{\qquad\qquad\! \st \quad 		\begin{bmatrix}
\Up \\ \Yp \\\Uf 
\end{bmatrix}
 g = \begin{bmatrix}
\uini \\ \yini \\ u 
\end{bmatrix} 
}.
\end{mini} 
We now follow the arguments from Section~\ref{subsec: multi-objective} to reduce the bi-level problem \eqref{eq:least-norm form} to a single-level multi-criteria problem.

As in \eqref{eq:OPT-BL-value}, the inner problem can be replaced by a constraint assuring that it achieves its minimum. 
{\tb Here, we add an orthogonality constraint to the constraints of the inner problem: 
\begin{equation*}
\begin{bmatrix}
\Up \\ \Yp \\\Uf 
\end{bmatrix}
 g = \begin{bmatrix}
\uini \\ \yini \\ u 
\end{bmatrix} 
\quad\text{and}\quad
0 = \| (I - \Pi) g \|_{p}
\end{equation*}}%
The {\tb orthogonality constraint $0 = \| (I - \Pi) g \|_{p}$} poses the inner optimality constraint as the distance to the subspace containing the minimizers of the inner problem.
Retaining all constraints, \eqref{eq:least-norm form} can then be formulated as the single-level problem
\begin{mini}%
{\tb u \in \mathcal U,y \in \mathcal Y,g}{\ctr\left( \left[\begin{smallmatrix} y - y_{r} \\ u - u_{r}   \end{smallmatrix} \right] \right)  }{\label{eq:ARX-SMO}}{}
\addConstraint{	
\begin{bmatrix}
\Up \\ \Yp \\\Uf \\ \Yf
\end{bmatrix}
 g = \begin{bmatrix}
\uini \\ \yini \\ u \\y
\end{bmatrix} 
}
\addConstraint{	
\| (I - \Pi) g \|_{p} = 0
}
\,.
\end{mini}
We now apply Proposition \ref{Proposition: LipschitzPenalty}, lift the distance constraint $\| (I - \Pi) g \|_{p} = 0$ to the objective, and recover  problem \eqref{eq:ARX-SMO-final} with $\lambda$ larger than the Lipschitz constant of $\ctr(\cdot)$. 

Hence, \eqref{eq:ARX-SMO-final} is equivalent to \eqref{eq:ARX-SMO} for $\lambda$ sufficiently large. Our final convex relaxation is to choose $\lambda$ small rather than large. Namely, from the view-point of the objective: it lowers the cost; or from the bi-level viewpoint: it turns the inner optimality constraint into a weaker sub-optimality constraint{\tb, i.e., we allow for solutions satisfying $\| (I - \Pi) g \|_{p} \geq  0$.}

{\tb Conclusion $(i)$ now follows since \eqref{eq:ARX-SMO-final} is convex; $(ii)$ follows since we have only enlarged the feasible set when passing from \eqref{eq:OPT-BL-ARX} to \eqref{eq:ARX-SMO-final}; and $(iii)$ follows due to the enlarged feasible set, since the costs of \eqref{eq:OPT-BL-ARX} and \eqref{eq:ARX-SMO}  coincide, and since \eqref{eq:ARX-SMO-final} is a relaxation of \eqref{eq:ARX-SMO} if $\lambda$ is not sufficiently large.}
\end{IEEEproof}

\begin{remark}[\tb Comments on Theorem~\ref{Theorem: SPC relaxation}]
\label{rem: comments on spc theorem}
%
First, we summarize the salient arguments to pass from indirect to direct data-driven control: we
relaxed problem \eqref{eq:OPT-BL-ARX} by dropping causality (block-triangularity) and LTI complexity (rank) specifications, replaced the least-square criterion \eqref{eq:pem} by the equivalent least-norm formulation \eqref{eq:least-norm form}, 
and lifted the problem from bi-level to multi-criteria, where the least-square objective {\tb induces} the regularization $\|(I-\Pi)g\|_{p}$. 
For equivalence to the least-square objective, the proof requires $\lambda$ larger than the (global) Lipschitz constant of $\ctr(\cdot)$, similar to  robustification-induced regularizations \cite{JC-JL-FD:19-CDC,JC-JL-FD:20}. If  $\ctr(\cdot)$ is only locally Lipschitz, e.g., in case of a quadratic cost, then choosing a finite (small) $\lambda$ is a relaxation that allows the predicted trajectory to not adhere to the least-square fit of the data. Though as we will see in Section~\ref{subsec: role of projection}, its effect is minor for $\lambda$ not overly small.

{\tb Second, continuing on the magnitude of $\lambda$: For exact data and under consistency assumptions, \eqref{eq:ARX-SMO-final} achieves the exact minimizer for  any $\lambda \geq 0$; see Fact~\ref{fact: consistency of projection}. When departing from these ideal assumptions, 
the least-square fit of the data is enforced only for $\lambda$ sufficiently large. Generally, $\lambda$ should be regarded as a tune-able hyper-parameter chosen by the designer to control how much the predicted trajectory should adhere to the data (versus the control objective) and to ultimately improve the realized performance. The proof of Theorem~\ref{Theorem: SPC relaxation} suggests a sufficiently large value, which is also confirmed by our later empirical findings (see e.g.  Figure~\ref{fig:projected_ell2_regularization_more_data}). 
}
 
{\tb Third}, the regularization based on the projector $\|(I-\Pi)g\|_{p}$ differs from the standard $p$-norm regularizers $h(g) = \|g\|_{p}$ \cite{xue2020data,JC-JL-FD:19-CDC,JC-JL-FD:20} (or squared 2-norms $\|g\|_{2}^{2}$ \cite{LH-JZ-JL-FD:20,berberich2020data}). Actually, it is this projection  which recovers 
 the least-square criterion \eqref{eq:pem}. {\tb In contrast, norm-based regularizers $\|g\|_{p}$ are not consistent and bias the optimal solution $(u^{\star},y^{\star})$; see Remark~\ref{rem: consistency}. 
This is undesirable from an identification perspective: the regularizer should induce a least-square fit of the data. While for small values of $\lambda$  both regularizers have a similar  effect, for sufficiently large $\lambda$  the {\em identification-induced regularizer} $\|(I-\Pi)g\|_{p}$ demonstrates a superior performance; see Figure~\ref{fig:projected_ell2_regularization_more_data} later.} 

{\tb Fourth}, our proof strategy reveals an entire class of regularizers. In fact, we can choose any $p$-norm $\|(I-\Pi)g\|_{p}$, use more general penalty functions such as the (squared) merit functions in \cite{ye1997exact}, or attack problem \eqref{eq:ARX-SMO} with other penalty or augmented Lagrangian methods. These degrees of freedom reflect the intuition that the Pareto-front of \eqref{eq:ARX-SMO-final} is invariant under certain (e.g., monotone) transformations of objectives such as taking squares; see \cite[Appendix A]{xu2010robust} for a formal reasoning. For our later simulations in Section~\ref{subsec: role of projection}, we choose the computationally attractive regularization $\|(I-\Pi)g\|_{2}^{2}$.

{\tb Fifth and finally,} our proof arguments are obviously ``qualitative'' crossing out {\tb rank and causality constraints similar to most SPC implementations} and using non-quantifiable ``sufficiently large'' reasoning. Hence, the convex relaxation \eqref{eq:ARX-SMO-final} of  \eqref{eq:OPT-BL-ARX} should not be expected to be tight. 
Nevertheless, the formulation \eqref{eq:ARX-SMO-final} (without projector) has proved itself {\tb  in many case studies} and often outperforms \eqref{eq:OPT-BL-ARX}, as testified in \cite{LH-JZ-JL-FD:01,LH-JC-JL-FD:19,PC-AF-SB-FD:20,EE-JC-PB-JL-FD:19,LH-JZ-JL-FD:20}. Section~\ref{subsec: role of projection} will compare the different formulations.
\oprocend
\end{remark}

\subsection{Bridging Towards Structured Low-Rank  Approximation}
\label{subsec: low-rank relaxation}

We now present an entirely non-parametric problem formulation, namely a version of subspace identification based on structured low-rank approximation \cite{markovsky2016}, and we relate the resulting bi-level problem to direct data-driven control \eqref{eq:OPT-DR}.

Given the model class $\lti_{m,\ell}^{q,n}$, we project the identification data $w_{d} \in \real^{qT}$ on $ \widehat\bv_{\Tini+\Tr} \in \lti_{m,\ell}^{q,n}$.  By Lemma \ref{lemma: fundamental lemma}, the latter set is characterized by all trajectories $\hat w \in \real^{q(\Tini+\Tr)}$ so that the associated Hankel matrix satisfies $\text{rank} \left(\Han_{\Tini+\Tr}(\hat w) \right) \leq m(\Tini+\Tr)+n$ for $(\Tini+\Tr) > \ell$. An implicit assumption is, of course, $\Td \gg \Tini+\Tr$: the identification data is much longer than the estimation plus control prediction horizons.

In presence of noise,  $\Han_{\Tini+\Tr}( w_{d}) $ will not have low rank and has to be approximated by a low-rank matrix in an identification step.
Thus, the identification problem \eqref{eq:ID} reads\,as%
\begin{mini}%
{\hat w_{d}}{\cid(\hat w_{d}-w_{d})}{\!\!\!\!\!\!\!\!\!\!\!\!\!\!\!\!\!\!\!\!\!\!\!\!\!\!\!\!\!\!\!\!\!\!\!\!\!\!\!\!\!\!\!\!\!\!\!\!\!\!}{\label{eq:low-rank}}
\addConstraint{
 \text{rank} \left(\Han_{\Tini+\Tr}(\hat w_{d}) \right) {\leq} m(\Tini+\Tr)+n.
}%
\end{mini}%
Problem \eqref{eq:low-rank} is to be read as low-rank approximation problem: given the identification data assorted in a Hankel matrix $\Han_{\Tini+\Tr}(w_{d})$, we seek the closest sequence $\hat w_{d}$ so that the Hankel matrix $\Han_{\Tini+\Tr}(\hat w_{d})$ has rank {no more than} $m(\Tini+\Tr)+n$.

Since $\hat w_{d} \in \widehat \bv_{T}$, we have  $\text{rank} \left(\Han_{\Tini+\Tr}(\hat w_{d}) \right) {\leq}  m(\Tini+\Tr)+n$. Since also  $\wini \in \widehat\bv_{\Tini}$ and $w \in \widehat\bv_{\Tr}$, we conclude
\begin{equation*}
\text{rank} \left(\left[\Han_{\Tini+\Tr}(\hat w_{d}) \,~\, \text{col}(\wini,w) \right]\right) {\leq} m(\Tini+\Tr)+n \,.
\end{equation*}
Assuming that $\text{rank} \left(\Han_{\Tini+\Tr}(\hat w_{d})\right) = m(\Tini+\Tr)+n$, which is generically the case, $\Han_{\Tini+\Tr}(\hat w_{d}) g = \text{col}(\wini,w) $ for some vector $g$. Hence, the bi-level problem \eqref{eq:OPT-BL} takes the form%
\begin{mini}
{w {\tb \in \mathcal W},g}{\!\!\!\!\!\ctr(w-w_{r}) }{\!\!\!\!\!\!\!\!\!\!\!\!\!\!\!\!\!\!\!\!\!\!\!\!\!\!\!\!\!\!\!\!\!\!\!\!\!\!\!\!\!\!\!\!\!\!\!\!\!\!\!\!\!\!}{\label{eq:low-rank-2}}
\addConstraint{\!\!\!\!\!
 \begin{bmatrix}
\wini\\ w\end{bmatrix} 
= \Han_{\Tini+\Tr}(\hat w_{d}^{\star})g 
}
\addConstraint{\!\!\!\!\!
\hat w_{d}^{\star} \in \argmin_{\hat w_{d}} \; \cid(\,\hat w_{d}-w_{d})
}
\addConstraint{\!\!\!\!\! \st\,  \textup{rank} (\Han_{\Tini+\Tr}(\,\hat w_{d})) \!=\! m(\Tini\!+\!\Tr)\!+\!n
}.
\end{mini}
%

\begin{theorem}[{\tb $\ell_{1}$-norm relaxation}] 
  \label{Theorem: Low-rank relaxation}
Consider the indirect data-driven control problem \eqref{eq:low-rank-2} and the direct data-driven control problem \eqref{eq:OPT-DR} for $h(g) = \|g\|_{1}$ and parameterized by $\lambda\geq0$. Let Assumptions~{\tb \ref{ass:ctr cost} and}~\ref{ass: cid} hold. For $\lambda$ sufficiently small, \eqref{eq:OPT-DR} is a convex relaxation of  \eqref{eq:low-rank-2}, that is,
\begin{enumerate}
\item[$(i)$] \eqref{eq:OPT-DR} is convex,
\item[$(ii)$] any feasible $(w,g)$ in \eqref{eq:low-rank-2} is also feasible for \eqref{eq:OPT-DR}, and 
\item[$(iii)$] the optimal value of \eqref{eq:OPT-DR} lower-bounds that of \eqref{eq:low-rank-2}.
\end{enumerate}
\end{theorem}

\begin{IEEEproof}
To prove the claim, one can resort to a proof strategy via the multi-criteria problem \eqref{eq:OPT-SMO}, as in the previous section. Instead, we present a more direct approach here.

We start by massaging the rank constraint in \eqref{eq:low-rank-2}. First, since $ \text{rank} \left(\Han_{\Tini+\Tr}(\hat w_{d})\right) = m(\Tini+\Tr)+n$, we may without loss of generality add the constraint $\|g\|_0 \leq n+ m(\Tini+\Tr)$ to the outer problem, where $\|g\|_0$ denotes the cardinality (number of nonzero entries) of $g$. Second, we perform a convex relaxation and drop the rank constraint. 
Third, another convex relaxation (popular in LASSO problems \cite{hastie2015statistical}) is to replace $\|g\|_0 \leq n+ m(\Tini+\Tr)$ by $\|g\|_{1} \leq \alpha$ for $\alpha>0$  sufficiently large. As a result of these three steps, \eqref{eq:low-rank-2} is relaxed to
\begin{mini}
{w{\tb \in \mathcal W},g}{\ctr(w-w_{r}) }{\!\!\!\!\!\!\!\!\!\!\!\!\!\!\!\!\!\!\!\!\!\!\!\!\!\!}{\label{eq:low-rank-3}}
\addConstraint{
 \begin{bmatrix}
\wini\\ w\end{bmatrix} 
= \Han_{\Tini+\Tr}(\hat w_{d}^{\star})g 
\,,\; \|g\|_1 \leq \alpha
}
\addConstraint{\!\where \;
\hat w_{d}^{\star} \in \argmin_{\hat w_{d}} \; \cid(\,\hat w_{d}-w_{d})
}.
\end{mini}
Observe that under Assumption~\ref{ass: cid} the inner problem admits a trivial solution:\, $\hat w_{d}^{\star}=w_{d}$. Thus, \eqref{eq:low-rank-3} reduces to 
\begin{mini}
{ w{\tb \in \mathcal W},g}{\ctr( w-w_{r}) }{\!\!\!\!\!\!\!\!\!\!\!\!\!\!\!\!\!\!\!\!\!\!\!\!\!\!\!\!\!\!\!\!\!\!\!\!\!\!\!\!\!\!\!\!\!\!\!\!\!\!\!\!}{\label{eq:low-rank-4}}
\addConstraint{
\begin{bmatrix}
\wini\\ w\end{bmatrix} 
= \Han_{\Tini+\Tr}( w_{d})g 
\,,\; \|g\|_1 \leq \alpha
}.
\end{mini}
Next, we lift the 1-norm constraint to the objective
\begin{mini}
{ w{\tb \in \mathcal W},g}{\ctr( w-w_{r}) + \lambda \cdot \|g\|_{1} }{\!\!\!\!\!\!\!\!\!\!\!\!\!\!\!\!\!\!\!\!\!\!\!\!\!\!}{\label{eq:low-rank-5}}
\addConstraint{
\begin{bmatrix}
\wini\\ w\end{bmatrix} 
= \Han_{\Tini+\Tr}( w_{d})g 
},
\end{mini}
where $\lambda\geq0$ is a scalar weight.
In particular, for each value of $\alpha$  in \eqref{eq:low-rank-4}, there is $\lambda\geq0$ so that the solution of \eqref{eq:low-rank-5} coincides with \eqref{eq:low-rank-4}, and vice versa. These equivalences are standard in $\ell_{1}$-regularized problems and follow from strong duality (applicable since $\ctr(\cdot)$ is convex and Slater's condition holds)  \cite{hastie2015statistical}.
 The precise value of $\lambda$ depends on the Lagrange multiplier of the constraint $\|g\|_1 \leq \alpha$ and thus  on the data. In either case, there is a selection of parameters so that both problems are equivalent, and choosing $\lambda$ sufficiently small is a relaxation. 
 
Thus, we arrived at the direct data-driven control \eqref{eq:OPT-DR} for $\lambda$ sufficiently small and $h(g) = \|g\|_{1}$.  
{\tb Conclusion $(i)$  follows due to convexity \eqref{eq:OPT-DR}; $(ii)$ follows since we have enlarged the feasible set passing from \eqref{eq:low-rank-2} to \eqref{eq:OPT-DR}; and $(iii)$ follows due to the enlarged feasible set, since the costs of \eqref{eq:low-rank-2} and \eqref{eq:low-rank-4}  coincide, and since \eqref{eq:OPT-DR} is a relaxation of \eqref{eq:low-rank-4} for $\lambda$ small.}
\end{IEEEproof}

In summary, to pass from indirect data-driven control \eqref{eq:low-rank-2} to direct data-driven control \eqref{eq:OPT-DR}, we performed a sequence of convex relaxations effectively replacing the rank constraint of the system identification by a $\ell_{1}$-norm regularizer. Hence, the 1-norm regularizer accounts for selecting the model complexity. Similar remarks as those following Theorem~\ref{Theorem: SPC relaxation} on tightness of the relaxation apply to Theorem~\ref{Theorem: Low-rank relaxation}, too; {\tb see Remark \ref{rem: comments on spc theorem}.}
%

\subsection{Hybrid relaxations}

Theorems \ref{Theorem: SPC relaxation} and \ref{Theorem: Low-rank relaxation} reveal the roles of the two regularizers: $\|g\|_{1}$ controls the {model} 
complexity, whereas $\|(I-\Pi)g\|_{2}$ accounts for least-square fitting the data. 
To blend the two, consider a hybrid formulation of \eqref{eq:ARX-SMO-final} and \eqref{eq:low-rank-2}
\begin{mini}
{w{\tb \in \mathcal W},g}{\!\!\!\!\!\ctr(w-w_{r})  \,+\, \lambda_{1} \cdot  \|(I-\Pi)g\|^{2}_{2}}{\!\!\!\!\!\!\!\!\!\!\!\!\!\!\!\!\!\!\!\!\!\!\!\!\!\!\!\!\!\!\!\!\!\!\!\!\!\!\!\!\!\!\!}{\label{eq:low-rank-6}}
\addConstraint{\!\!\!\!\!
 \begin{bmatrix}
\wini\\ w\end{bmatrix} 
= \Han_{\Tini+\Tr}(\hat w_{d}^{\star})g 
}
\addConstraint{\!\!\!\!\!
\hat w_{d}^{\star} \in \argmin_{\hat w_{d}} \; \cid(\,\hat w_{d}-w_{d})
}
\addConstraint{\!\!\!\!\! \st\,  \textup{rank} (\Han_{\Tini+\Tr}(\,\hat w_{d})) \!=\! m(\Tini\!+\!\Tr)\!+\!n
},
\end{mini}
where $\lambda_{1} \geq 0$. Observe that this formulation is consistent:
\begin{fact}
Under Assumptions \ref{ass:L+n pe} with $L$ replaced by $\Tini+L$, \ref{ass:ctr cost}, \ref{ass: Bp = Br}, and \ref{ass: Bp=Bid}, for any $\lambda \geq 0$ the minimum of \eqref{eq:low-rank-6} and achieved for  $w^{\star} = w_{r}$ and $\|(I-\Pi)g^{\star}\|^{2}_{2}=0$.
\end{fact}
{\tb The arguments in the previous section then lead us to} 
\begin{mini}
{ w {\tb\in\mathcal W},g}{\!\!\!\ctr( w-w_{r})  \,+\, \lambda_{1} \cdot  \|(I-\Pi)g\|^{2}_{2}  \,+\, \lambda_{2} \cdot  \|g\|_{1}}{\!\!\!\!\!}{\label{eq:low-rank-7}}
\addConstraint{\!\!\!
\begin{bmatrix}
\wini\\ w\end{bmatrix} 
= \Han_{\Tini+\Tr}( w_{d})g 
},
\end{mini}
where $\lambda_{2} \geq 0$. We will validate the performance of the hybrid regularizer in Section~\ref{subsec: role of projection} below; see specifically Figure~\ref{fig:hybrid_regularization_more_data}.

\subsection{Possible pitfalls of relaxations}
\label{subsec: Pitfall of relaxations}

Note that the two convex relaxation results in Theorems \ref{Theorem: SPC relaxation} and \ref{Theorem: Low-rank relaxation} are {\em trivially} true in the limit when $\lambda = 0$. In fact, even the abstract multi-criteria formulation \eqref{eq:OPT-SMO} can be related to a  relaxation of the abstract bi-level problem \eqref{eq:OPT-BL}  in the limit $\gamma = 0$. Namely, for $\gamma = 0$,  \eqref{eq:OPT-SMO}   reduces to%
\begin{mini}%
{ w,\hat w_{d},\widehat\bv}{\ctr( w-w_{r}) }{\!\!\!\!\!\!\!\!\!\!\!\!\!\!\!\!\!\!\!\!\!\!\!\!\!\!\!\!\!\!\!\!\!\!\!\!\!\!\!\!\!\!\!\!\!\!\!\!\!\!\!\!\!\!\!\!\!\!\!\!\!\!\!\!\!\!\!\!\!\!\!\!\!\!\!\!\!\!\!\!\!\!\!\!\!\!}{\label{eq:OPT-SMO-R}}
\addConstraint{
{\wini \wedge  w \in \widehat\bv_{\Tini+\Tr}}\,,\; \hat w_{d} \in \widehat\bv_{\Td}\,,\; \widehat\bv \in  \lti_{m,\ell}^{q,n}
\,.
}
\end{mini}%
{\tb The variable $\hat w_{d}$ and the constraint $\hat w_{d} \in \widehat\bv_{\Td}$ can be removed, and \eqref{eq:OPT-SMO-R} amounts to matching the model $\widehat\bv$ to the reference $w_{r}$. The next result is followed by a discussion on regularizers:}

\begin{corollary}\label{Corollary: trivial relaxation}
Consider the indirect data-driven control \eqref{eq:OPT-BL} and multi-criteria problem \eqref{eq:OPT-SMO-R} in the limit $\gamma = 0$, {\tb and let Assumption~\ref{ass:ctr cost} hold.} Then problem \eqref{eq:OPT-SMO-R} is a  relaxation of problem \eqref{eq:OPT-BL}, that is,
\begin{enumerate}
\item[$(i)$] any feasible $( w,\hat w_{d},\hat \bv)$ in \eqref{eq:OPT-BL} is also feasible for \eqref{eq:OPT-SMO-R},  
\item[$(ii)$] and the optimal value of \eqref{eq:OPT-SMO-R} lower-bounds that of \eqref{eq:OPT-BL}.
\end{enumerate}
\end{corollary}

\begin{IEEEproof}
Consider the equivalent formulation \eqref{eq:OPT-BL-value} of \eqref{eq:OPT-BL}, and note that \eqref{eq:OPT-SMO-R} equals \eqref{eq:OPT-BL-value} when the inner optimality constraint $\cid(\hat w_{d}-w_{d}) - \varphi = 0$ is dropped. {\tb The conclusions now follow analogously as in Theorems \ref{Theorem: SPC relaxation} and \ref{Theorem: Low-rank relaxation}.}
\end{IEEEproof}

Analogous corollaries can be stated for Theorems \ref{Theorem: SPC relaxation} and \ref{Theorem: Low-rank relaxation} for $\lambda = 0$. 
Given such results, one may wonder whether Theorems \ref{Theorem: SPC relaxation} and \ref{Theorem: Low-rank relaxation} are vacuous since they are trivially true for $\lambda = 0$. We offer several answers. 
First, the limit $\lambda =0$ clearly leads to a better solution $w^{\star}$ (i.e., a lower surrogate tracking error) for the {\em open-loop} optimal control problem. However, this solution merely matches the reference $w_{r}$ and does not adhere to the identification data $w_{d}$ in the sense of meeting any fitting criterion. Hence, the optimal solution $w^{\star}$ may not be a trajectory of the true system behavior, and the actual {\em realized} control performance can be arbitrarily poor. Obviously, such a situation is not desirable, and one may want to regularize with a small but non-zero $\lambda$ -- an observation consistent with \cite{JC-JL-FD:18,JC-JL-FD:19-CDC,JC-JL-FD:20,xue2020data,LH-JZ-JL-FD:20,LH-JZ-JL-FD:01} albeit derived from a different perspective.
Second, Theorems \ref{Theorem: SPC relaxation} and \ref{Theorem: Low-rank relaxation} require $\lambda$ to be sufficiently small, but not  zero. According to the proofs, depending on Lipschitz constants and multipliers of the respective  problems, there is a smallest value for $\lambda$ so that the behavior $\widehat\bv$ matches (in the $\cid(\cdot)$ fitting criterion) the plant behavior $\bv^{P}$. 
In \cite{JC-JL-FD:18,JC-JL-FD:19-CDC,JC-JL-FD:20,xue2020data,LH-JZ-JL-FD:20,LH-JZ-JL-FD:01} the coefficient $\lambda$ relates to a desired robustness level. In either case, $\lambda$ can hardly be quantified a priori and without cross-validation; {\tb see also Remark~\ref{rem: comments on spc theorem}.}

We follow up on this set of questions in the next section.

\section{Numerical Analysis and Comparisons}
\label{subsec: numerical analysis}

We now numerically investigate the effect of the hyper-parameter $\lambda$, confirm the superiority of the regularizer $h(g) = \|(I-\Pi)g\|_{2}^{2}$, and compare direct and indirect approaches.

\subsection{Choice of Regularization Parameter}
\label{subsec: choice of regularization}

We first study the parameter $\lambda$ regularizing direct data-driven control \eqref{eq:OPT-DR}. 
Consider the benchmark single-input, single-output, 5th order, linear time-invariant system \cite{ddctr-benchmark}. 
Denoting the $t$-th element of the concatenated input and output by ${w}(t)=({u}(t),{y}(t))$, the control cost was chosen as $c_\textup{ctrl}({w}-w_r) = ({w}-w_r)^{\top}W({w}-w_r)$ with reference $w_r(t) = (u_r(t),y_r(t)) = (0,\sin(2\pi t/(L-1)))$ for $t\in\{0,1,\dots,L-1\}$, prediction horizon $L=20$, $W=I_L\otimes\textup{diag}(0.01,2000)$, where $I_L$ is the $L \times L$ identity, and $\otimes$ denotes the Kronecker product. {\tb In this entire section, we disregard constraints, i.e., $\mathcal W \equiv \real^{q\Tr}$.} We used a 1-norm regularizer $h(g) = \|g\|_{1}$ in~\eqref{eq:OPT-DR} and a prefix-trajectory of length $\Tini=5$ (see Section~\ref{subsec: estimation}).

\begin{figure}[b]
    \centering 
  \includegraphics[width=\columnwidth]{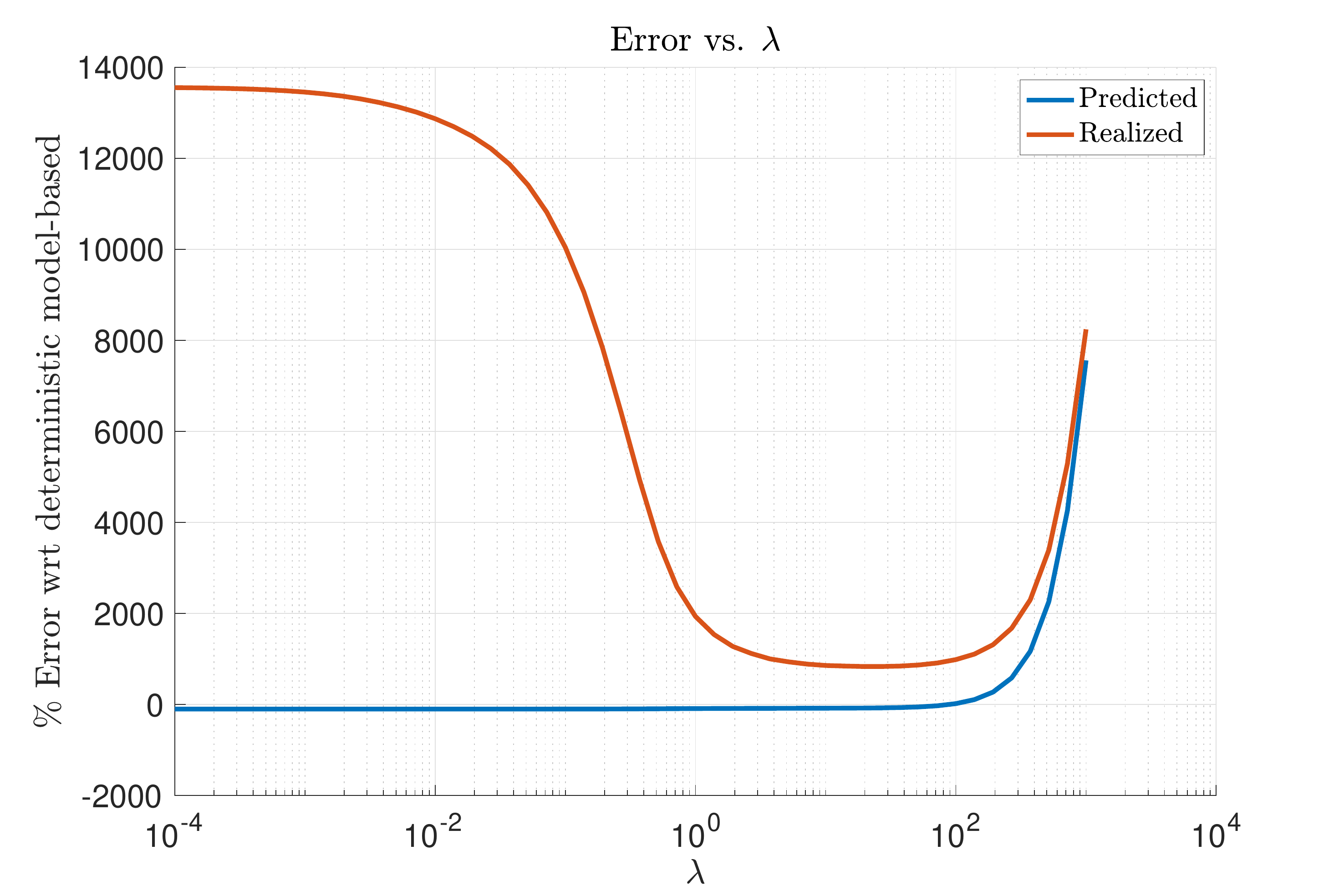}
  \caption{Predicted and realized errors {\tb (relative to the ground-truth optimal performance and averaged over 100 data sets) with 1-norm regularizer $\lambda \|g\|_{1}$.}}
  \label{fig:error_vs_epsilon}
\end{figure}

We collected one noise-free input/output time series of length {\tb$T=250$} by applying a random Gaussian input. From this noise-free data set, 100 independent noisy data sets were constructed by adding Gaussian noise with a noise-to-signal ratio of 5\%. For each data set and each value of $\lambda\in(0,10^3)$, optimal control inputs were computed from~\eqref{eq:OPT-DR}. We define the {\em predicted} error as $c_{\textup{ctrl}}(w^{\star}-w_r)$, where $w^{\star}$ is an optimizer of \eqref{eq:OPT-DR}. We define the {\em realized} error as $c_{\textup{ctrl}}(w_{\textup{true}}-w_r)$, where $w_{\textup{true}}$ is the realized trajectory of the system after applying the computed optimal inputs. The predicted and realized errors were converted to a percentage increase in error with respect to the ground-truth optimal performance (i.e., if the deterministic system was exactly known), and were averaged over the 100 independent data sets. The results are plotted in Figure~\ref{fig:error_vs_epsilon}. 

It is apparent that choosing $\lambda$ too small leads to an optimistic predicted error but very poor realized performance. Furthermore, the performance is poor for  large values of $\lambda$ indicating that the regularization parameter should be chosen carefully (though a wide range delivers equally good results). These observations are consistent with those in~\cite{JC-JL-FD:18,xue2020data,JC-JL-FD:19-CDC,JC-JL-FD:20,LH-JZ-JL-FD:20,LH-JZ-JL-FD:01} and the hypotheses discussed at the end of Section~\ref{subsec: Pitfall of relaxations}.
%
%
%

\subsection{Role of Projection in Two-Norm Regularization}
\label{subsec: role of projection}

\begin{figure}[h]
    \centering 
  \includegraphics[width=\columnwidth]{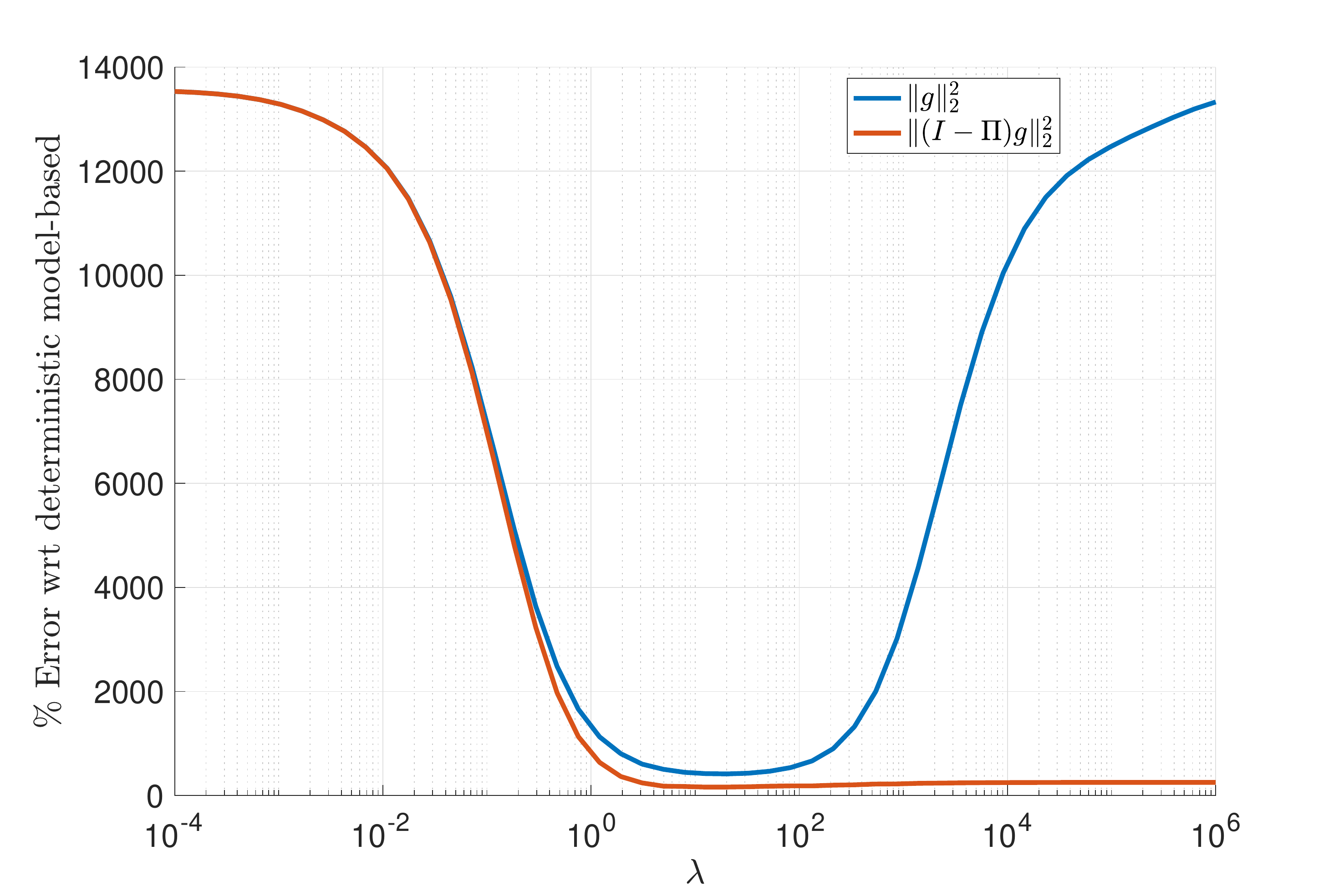}
  \caption{Comparison of the realized performance {\tb (relative to the ground-truth optimal performance and averaged over 100 data sets)} for the two-norm $\|g\|_{2}^{2}$ and identification-induced regularization $\|(I-\Pi)g\|_{2}^{2}$ as function of $\lambda$.}
  \label{fig:projected_ell2_regularization_more_data}
\end{figure}

Theorem~\ref{Theorem: SPC relaxation} suggests that the identification-induced regularizer $h(g) = \|(I-\Pi)g\|_{2}^{2}$ is superior to a two-norm regularizer $h(g) = \|g\|_{2}^{2}$ if one is interested in consistency and the predicted trajectory adhering to a least-square fit of the data. 
To test this hypothesis, we consider the same case study from Section~\ref{subsec: choice of regularization} 
and report the averaged cost in Figure~\ref{fig:projected_ell2_regularization_more_data}. 
%
%

Both regularizers perform similarly for small $\lambda$, but the identification-induced regularizer shows a superior and surprisingly constant performance for sufficiently large  $\lambda$. By the proof of Theorem~\ref{Theorem: SPC relaxation}, for $\lambda$ sufficiently large, the direct and indirect problems \eqref{eq:ARX-SMO-final} and \eqref{eq:OPT-BL-ARX} are equivalent -- up to causality and complexity constraints. Thus, a sufficiently large $\lambda$ forces the least-square fit \eqref{eq:pem} and results in excellent performance independent of the specific value of $\lambda$. While there is a small window where the two-norm excels, the identification-induced regularizer shows overall much more robust performance.

\begin{figure}[tb]
    \centering 
  \includegraphics[width=\columnwidth]{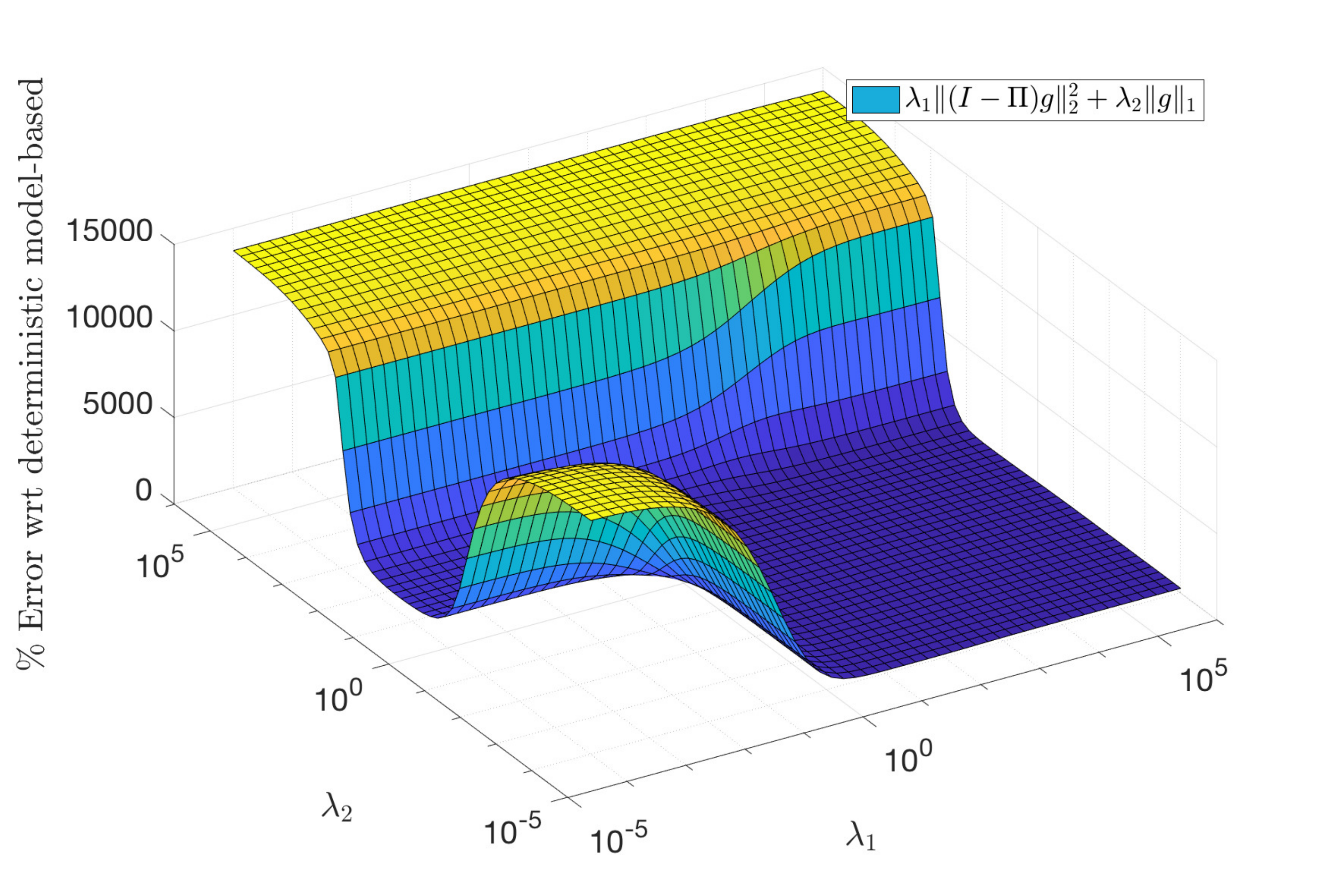}
  \caption{Realized error {\tb (relative to the ground-truth optimal performance and averaged over 100 data sets)} for a hybrid regularizer $\lambda_{1}   \|(I-\Pi)g\|_{2}^{2}  + \lambda_{2}   \|g\|_{1}$}
  \label{fig:hybrid_regularization_more_data}
\end{figure}
Next we study the merits of hybrid regularization \eqref{eq:low-rank-7}. For the same case study Figure~\ref{fig:hybrid_regularization_more_data} shows the averaged realized performance plotted over the regularization parameters. The $\{\lambda_{1}=0\}$ and $\{\lambda_{2}=0\}$ slices recover Figures \ref{fig:error_vs_epsilon} and \ref{fig:projected_ell2_regularization_more_data}. As before, the regularizer $ \|(I-\Pi)g\|_{2}^{2}$ is more robust though a hybrid regularizer yields a minor albeit robust improvement. {\tb A closer examination of the data underlying Figure~\ref{fig:hybrid_regularization_more_data} reveals that a hybrid regularization can improve up to 15\% over the best results achievable with the regularizer  $ \|(I-\Pi)g\|_{2}^{2}$  only.}

{\tb
\subsection{Effect of data length}

We continue with the same case study and discuss the effect of data-length on direct and indirect methods. 
For the direct method, we used the identification-induced regularizer $h(g) = \|(I-\Pi)g\|_2^2$ with sufficiently large weight $\lambda = 10000$, as indicated in Figure \ref{fig:projected_ell2_regularization_more_data}. For the indirect method, the inner system identification problem~\eqref{eq:ID} is solved using the subspace approach N4SID~\cite{van1994n4sid} with prefix horizon $\Tini=5$, prediction horizon $L=20$, and (correct) model-order selection $n=5$. 

For our case study Lemma~\ref{lemma: fundamental lemma} demands at least $T = 59$ data points. Figure~\ref{fig:performance_vs_numdata} below shows the beneficial effects of including more data on the realized {\em median} performance of the direct and indirect methods. The main findings are as follows: First,  both methods are asymptotically consistent. Second, the indirect method is superior in the low data regime echoing that models are compressed and de-noised  representations, see Remark~\ref{rem: models vs. data}. Third and finally, when an incorrect model-order $n=6$ is selected for the indirect method (resulting in an over-parameterization and thus a bias), then consistency is lost, and the direct method is superior. This effect is even more pronounced when studying the average (as opposed to the median) error due to several outliers of the indirect method.

This third point hints at a bias-variance trade-off between the direct and indirect methods, which will be studied below.

\begin{figure}[h]
    \centering 
  \includegraphics[width=\columnwidth]{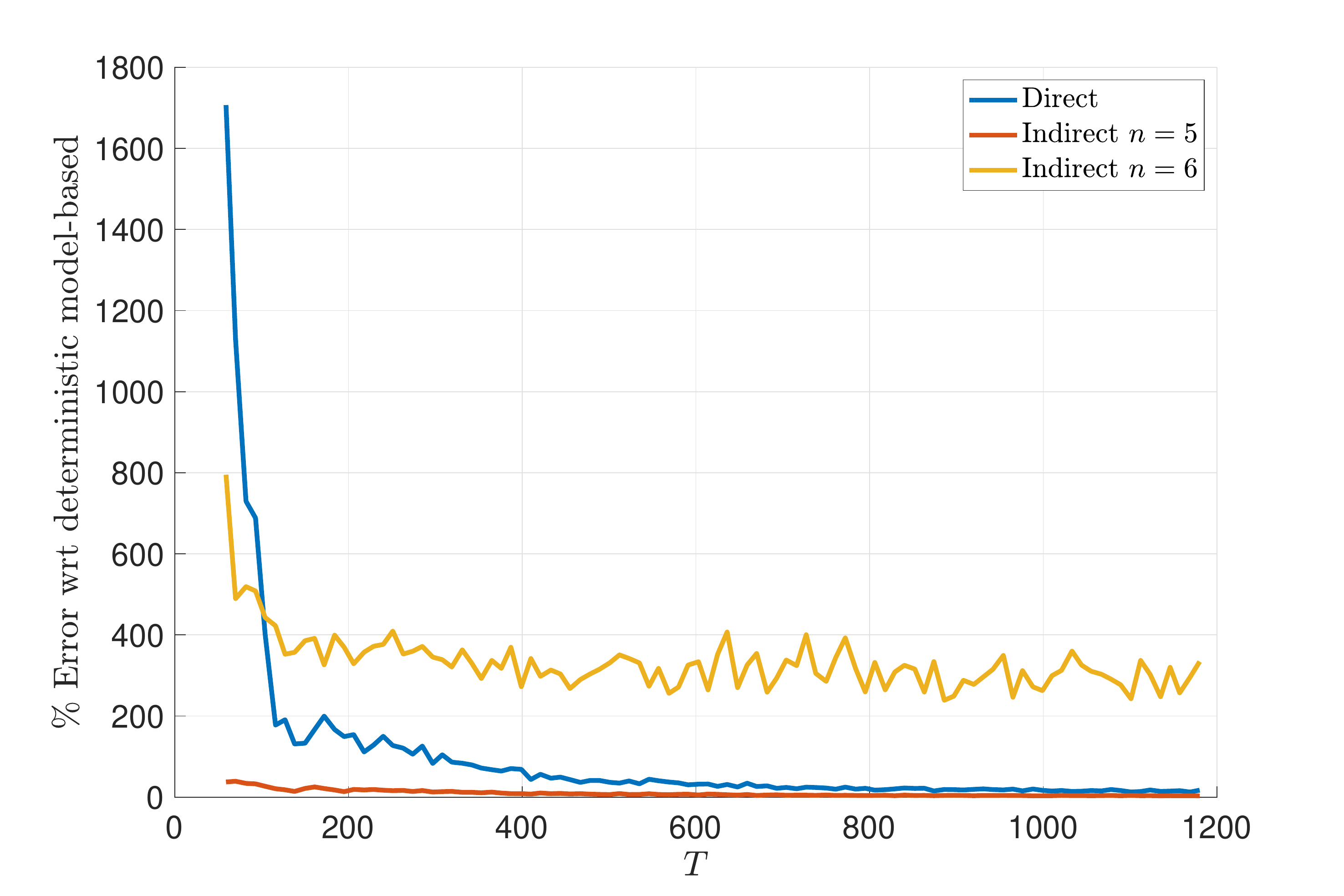}
  \caption{Realized median error (over 100 data sets) for the direct and indirect (with different model order selections) methods for varying amount of data.}
  \label{fig:performance_vs_numdata}
\end{figure}

}

\subsection{Comparison and Bias-Variance Hypotheses}

We now compare the direct and indirect approaches through two case studies. The first study evaluates the performance of both methods on the basis of ``variance'' error, i.e., on a linear system with noisy measurements. The second study evaluates the performance on the basis of ``bias'' error, i.e., on a nonlinear system with noise-free measurements. 

We expect the direct method to perform better on the nonlinear system since the indirect method erroneously selects a linear model class thus leading to a larger ``bias'' error. On the other hand, we expect the indirect method to perform better on the linear system with noisy outputs since the identification step filters  noise thus  leading to a lower ``variance'' error.

\subsection*{Comparison: Stochastic Linear System}

Consider the same case study as in the  Section~\ref{subsec: choice of regularization}, i.e., same LTI system, cost, and reference.
We collected data for varying levels of noise-to-signal ratio, i.e., we considered measurements that were affected by Gaussian noise with noise-to-signal ratio in the set $\{0\%,1\%,\dots,15\%\}$. For each noise-to-signal ratio, {\tb $T=250$} input/output data samples were collected by applying a random Gaussian input. This data was then used for both the direct and indirect methods.

For the indirect method, the inner system identification problem~\eqref{eq:ID} is {\tb again} solved using 
N4SID~\cite{van1994n4sid} with prefix horizon $\Tini=5$ and prediction horizon $L=20$. Equipped with a  (correct) 5th-order identified model, optimal control inputs are computed by solving~\eqref{eq:OPT-CE}. The indirect method was compared to the direct method~\eqref{eq:OPT-DR}, with $h(g) = \|g\|_{1}$, $\Tini=5$, and $\lambda=27$. The hyper-parameters of both methods were kept constant for all simulations below and chosen to give good realized control  performance for all noise-to-signal ratios.

\begin{figure}[h]
    \centering 
  \includegraphics[width=\columnwidth]{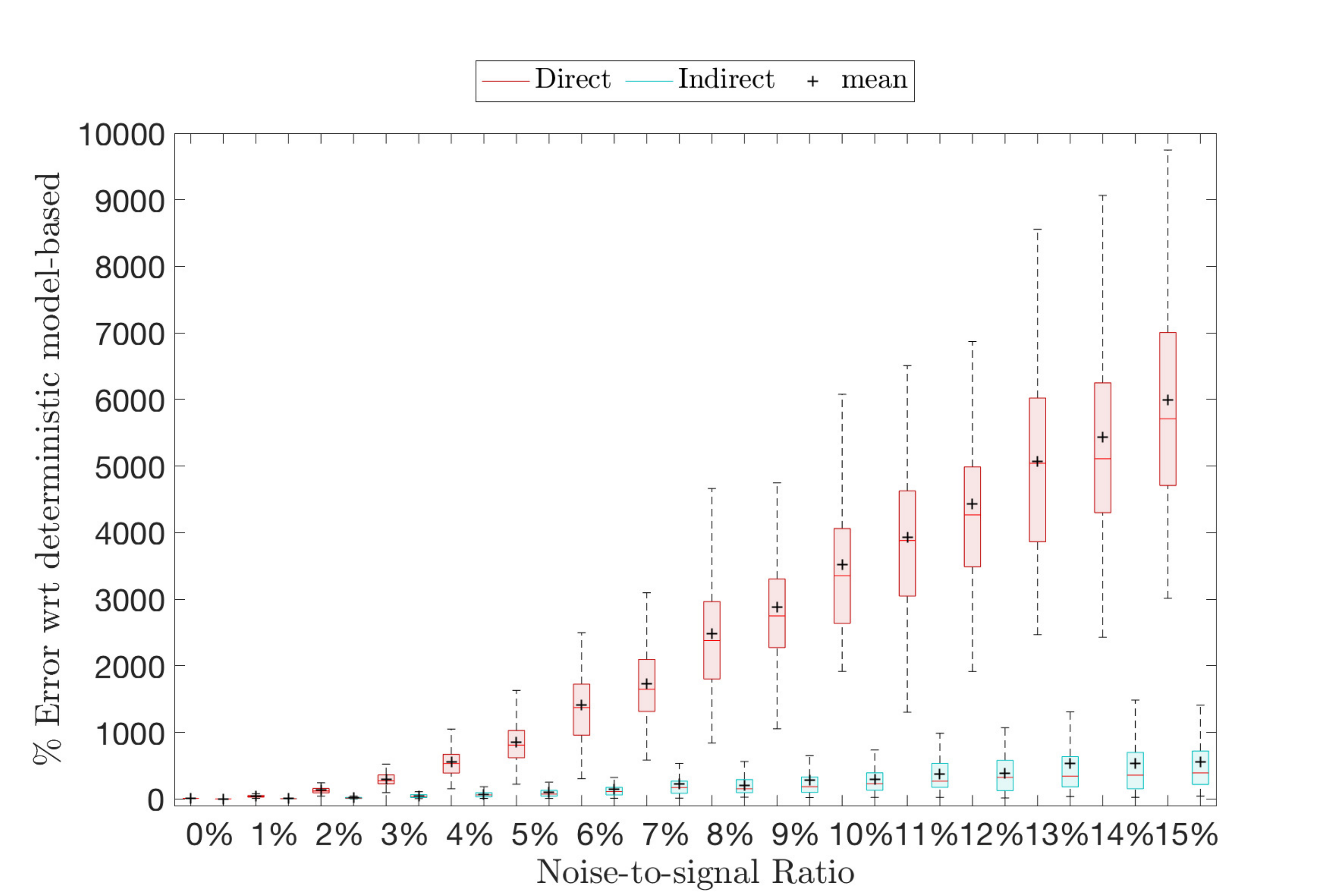}
  \caption{Comparison of direct and indirect methods for varying noise.}
  \label{fig:boxplot_noisy}
\end{figure}

For both methods we recorded the realized performance after applying the open-loop inputs and converted it to a percentage error with respect to the best possible performance (i.e., if the deterministic model was exactly known). For each noise-to-signal ratio, 100 simulations were conducted with different random data sets. The results are displayed in the box plot in Figure~\ref{fig:boxplot_noisy} and show that both methods perform well for low levels of noise (up to approximately $2\%$ noise-to-signal ratio). As the data becomes noisier, the performance of the direct method degrades significantly, while the performance of the indirect method remains relatively constant.  We remark that a slightly better albeit qualitatively similar result is obtained with the regularizer  $ \|(I-\Pi)g\|_{2}^{2}$.

We attribute these observations to the fact that identification de-noises the data. These results confirm our hypothesis that the indirect method is superior in terms of ``variance'' error.

%
%
%

\subsection*{Comparison: Deterministic Nonlinear System}

We now consider the scenario where the direct and indirect methods are subject to a ``bias'' error, but not a ``variance'' error. Consider the discrete-time nonlinear Lotka-Volterra dynamics  considered for direct data-driven control in~\cite{kaiser2018sparse}
\[
\begin{aligned}
x(t_{k+1}) 
&= f_{\textup{nonlinear}}(x(t_k),u(t_k))
\\
&=\left[\begin{smallmatrix}
 x_1(t_k) + \Delta t(ax_1(t_k) - bx_1(t_k)x_2(t_k)) \\
 x_2(t_k) + \Delta t(dx_1(t_k)x_2(t_k) -cx_2(t_k) + u(t_k))
\end{smallmatrix}\right]\,,
\end{aligned}
\]
where $t_{k+1} - t_k = \Delta t = 0.01$, $a=c=0.5$, $b=0.025$, $d=0.005$,  and $x(t_k) = \begin{bmatrix} x_1(t_k) & x_2(t_k)\end{bmatrix}^{\top}$. Here, $(x_1(t_k),x_2(t_k))$ denote prey and predator populations, and $u(t_k)$ is the input. A linearization about the equilibrium  $(\bar{u},\bar{x}_1,\bar{x}_2)=(0,c/d,a/b)$ 
yields the affine {\tb linear} system
\[
\begin{aligned}
&x(t_{k+1}) 
= f_{\textup{linear}}(x(t_k),u(t_k),\bar{x}_1,\bar{x}_2)
\\
&= \left[\begin{smallmatrix}
 x_1(t_k) + \Delta t\left((a-b\bar{x}_2)(x_1(t_k)-\bar{x}_1) - b\bar{x}_1(x_2(t_k)-\bar{x}_2)\right) \\
 x_2(t_k) + \Delta t\left(d\bar{x}_2(x_1(t_k)-\bar{x}_1) +(d\bar{x}_1 -c)(x_2(t_k)-\bar{x}_2) + u(t_k)\right)
\end{smallmatrix}\right]\,.
\end{aligned}
\]
{\tb We expect direct data-driven control \eqref{eq:OPT-DR} to perform well on such a nonlinear system for two reasons: $(i)$ nonlinear systems can be well approximated by LTI systems of sufficiently high complexity; and $(ii)$ the direct method \eqref{eq:OPT-DR} does not specify the LTI system complexity (e.g., by enforcing rank constraints).}

We compare the direct and indirect methods for varying degree of nonlinearity by interpolating between $f_{\textup{nonlinear}}$ and $f_{\textup{linear}}$, i.e., we study the interpolated system
\begin{equation}
\begin{aligned}
x(t_{k+1}) &= \epsilon \cdot f_{\textup{linear}}(x(t_k),u(t_k),\bar{x}_1,\bar{x}_2)\\
&\quad + (1-\epsilon) \cdot f_{\textup{nonlinear}}(x(t_k),u(t_k))
\end{aligned}
\label{eq:interpolated_dynamics}
\end{equation}
for $\epsilon\in [0,1]$. For $\epsilon=1$ (resp. $\epsilon=0$), the dynamics are purely affine (resp. nonlinear). For each $\epsilon\in\{0,0.1,\dots,1\}$, $T=2415$ data points were collected by applying a noisy sinusoidal input $u(t_k) = 2(\sin(t_k)+\sin(0.1t_k)))^2 + v(t_k)$ with $v(t_k)$ sampled from a Gaussian random variable. Full state measurement was assumed. The data collection was repeated for 100 different initial conditions. For each degree of nonlinearity $\epsilon\in\{0,0.1,\dots,1\}$ and each initial condition, the data was used to compute optimal open-loop control inputs using direct and indirect methods. The control cost was chosen as $c_\textup{ctrl}({w}-w_r) = \| {w}-w_r\|_{2}^{2}$ with equilibrium reference $w_r = (0,100,20,\dots,0,100,20)$, $L=600$, and $w=(u,x)$.

For the indirect method, the inner system identification optimization problem given by~\eqref{eq:ID} is solved using the subspace approach N4SID~\cite{van1994n4sid} with initial condition horizon $\Tini=4$, and prediction horizon $L=600$. A model order of 4 was chosen, as it produced the best performance as measured by the realized control cost. 
Optimal control inputs were then computed by solving~\eqref{eq:OPT-CE}. 
For comparison, we chose the direct method~\eqref{eq:OPT-DR} with $h(g) = \|g\|_{1}$, $\Tini=4$, and $\lambda=8000$. 
The performance was measured with the realized control cost after applying the open-loop inputs to system~\eqref{eq:interpolated_dynamics}. As before, the hyper-parameters of both direct/indirect methods were judicially chosen and kept constant for all simulations.%
\begin{figure}[tb]
    \centering 
\includegraphics[width=\columnwidth,trim=2cm 0cm 2cm 0cm,clip=true]{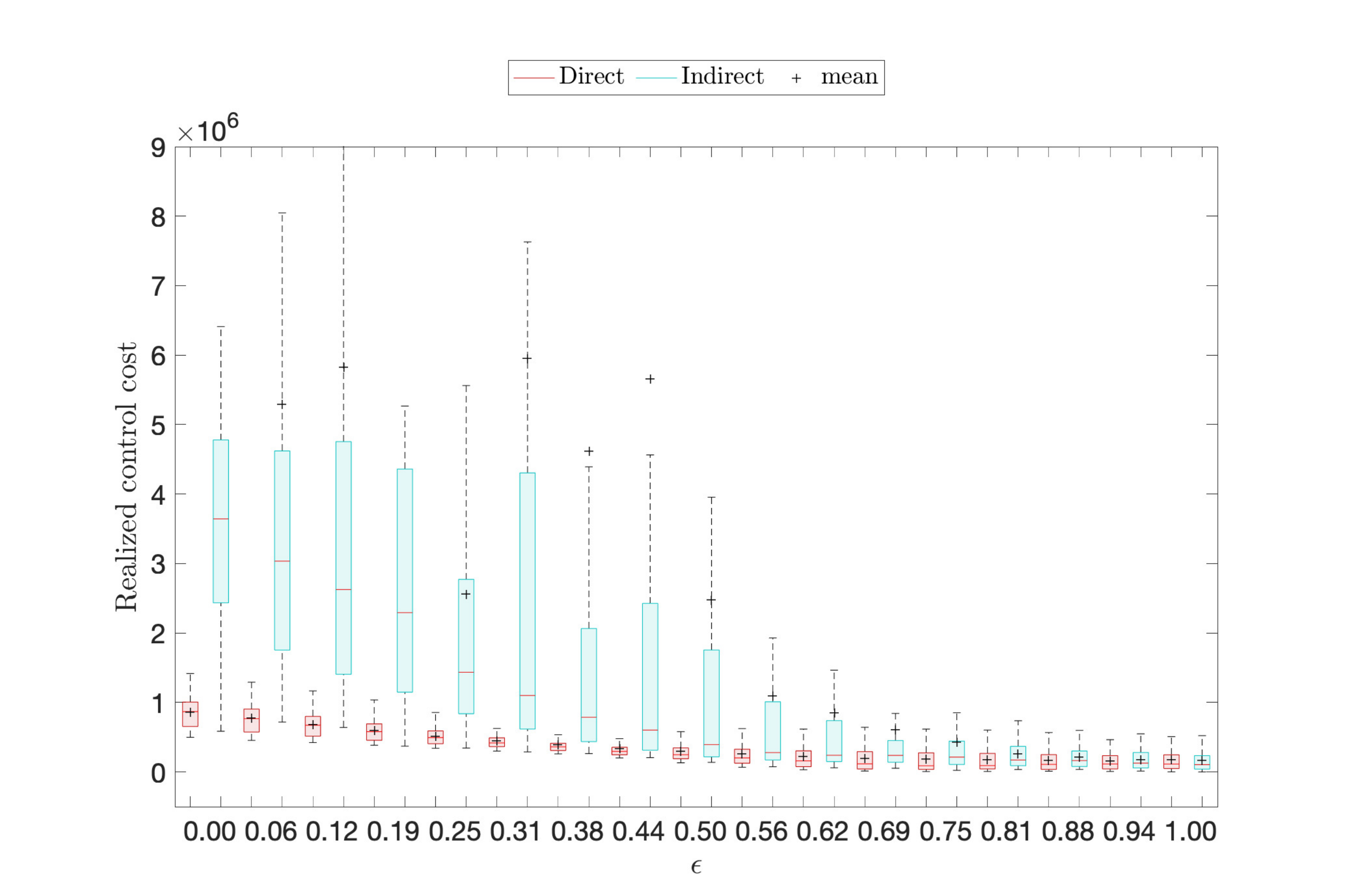}
  \caption{Comparison of direct and indirect methods for varying nonlinearity.}
  \label{fig:boxplot_nonlinear}
\end{figure}
%

The results displayed in Figure~\ref{fig:boxplot_nonlinear} show that both methods perform well for low levels of nonlinearity: $\epsilon \in [0.7, 1]$. As the system becomes increasingly nonlinear, the performance of the indirect method degrades significantly, while the performance of the direct method remains relatively constant. We attribute this observation to the fact that the indirect method incurs a ``bias'' error from selecting a linear model class and applying certainty-equivalence control, while the direct method uses data from the nonlinear system without bias. 
{\tb These findings confirm our earlier bias-variance observations from Figure~\ref{fig:performance_vs_numdata}.} 


\section{Discussion and Conclusions}
\label{sec: conclusions}

We studied the relationship between indirect and direct data-driven control formulated as bi-level (first-identify, then control) and single-level regularized (based on the Fundamental Lemma) optimization problems, respectively.  An intermediate multi-criteria problem allowed us to efficiently transition between both formulations. We concluded that the regularized direct approach can be viewed as a convex relaxation of the indirect approach, where the choice regularizer depended on the problem formulation and accounted for an implicit identification step. We also discovered a novel regularizer that is consistent and accounts for least-square identification. 

Our results suggested the use of the indirect method in case of ``variance'' errors and the use of the direct method in presence of ``bias'' errors (e.g., a nonlinear system {\tb or when selecting a wrong model order}). These insights {\tb echo the bias-variance trade-offs previously encountered for direct and indirect methods in \cite{campestrini2017data,krishnan2021direct}, and they}
 shed some partial light on the remarkable {\tb empirical performance of (direct)} data-enabled predictive control applied to nonlinear systems.  

As a  limitation, our results  concern only the open-loop predictive control problem, though we ultimately care about the realized performance, especially in a receding horizon closed-loop implementation. {\tb Some preliminary results on the realized performance of regularized control formulations were obtained in \cite{LH-JZ-JL-FD:01} through the lens of robust optimization, but the topic remains largely open.} Moreover, we believe that the proposed multi-criteria data-driven control formulation is important in its own right and may deliver excellent performance if one were to find a convex formulation and appropriate trade-off parameter. Both of these are formidable tasks for future work. 

Finally, we believe that our approach is also applicable to other identification and control formulations and may deliver interesting and novel direct data-driven control formulations.


\section*{Acknowledgements}

The authors acknowledge their colleagues at ETH Z\"urich, in particular Miguel Picallo Cruz, for fruitful discussions.


\appendix

Consider the mathematical program {\tb(MP)} 
\begin{mini}
{x  \in C}{f(x)}{\label{eq:MP}}{\!\!\!\!\!\!\!\!\!\!\!\!\!\!\!\!\!\!\!\!\!\!\!\!\!\!\!\!\!\!\!\!\!\!\!\!\boldsymbol{MP}:}
\addConstraint{
g(x) \leq 0
\,,\,
h(x) = 0
\,,
}
\end{mini}
where $C \subset \real^{n}$ is closed, and $f,g,h$ are lower semicontinuous maps from $\real^{n}$ to $\real$, $\real^{m}$, and $\real$. Consider the perturbation
\begin{mini}
{x  \in C}{f(x)}{\label{eq:MPeps}}{\!\!\!\!\!\!\!\!\!\!\!\!\!\!\!\!\!\!\!\!\!\!\!\!\!\!\!\!\!\!\!\!\!\!\!\!\boldsymbol{MP}_{\epsilon}:}
\addConstraint{
g(x) \leq 0
\,,\,
h(x) = \epsilon
}
\end{mini}
for $\epsilon \in \real$. We recall the definition of partial calmness \cite{ye1997exact}: 

\begin{definition}
 Let $x^{\star}$ solve $\boldsymbol{MP}$, and let $\mathbb B_{n}$ denote the open unit ball in $\mathbb R^{n}$. Then $\boldsymbol{MP}$ is said to be {\em partially calm}  at $x^{\star}$ provided that there are $\mu>0$ and $\delta>0$ such that, for all $\epsilon \in \delta \mathbb B_{1}$ and all $x \in x^{\star}+\delta \mathbb B_{n}$ feasible for $\boldsymbol{MP}_{\epsilon}$, one has 
\begin{equation}
f(x) + \mu | h(x) | \geq f(x^{\star})
\label{eq: calmness}
\,.
\end{equation}
\end{definition}
Partial calmness is equivalent to {\em exact penalization}. In particular, consider for $\mu \geq 0$ the penalized mathematical program 
\begin{mini}
{x  \in C}{f(x) \,+\, \mu \cdot |h(x)|}{\label{eq:PMP}}{\boldsymbol{PMP}_{\mu}:}
\addConstraint{
g(x) \leq 0
\,.
}
\end{mini}
We summarize \cite[Proposition 3.3]{ye1995},\cite[Proposition 2.2]{ye1997exact}:

\begin{proposition}\label{Proposition: partial calmness and exact penalty}
Assume that $f$ is continuous, and let $x^{\star}$ be a local minimizer of $\boldsymbol{MP}$.
Then $\boldsymbol{MP}$ is partially calm at $x^{\star}$ if and only if there is $\mu^{\star} >0$ so that $x^{\star}$ is a local minimizer of $\boldsymbol{PMP}_{\mu}$ for all $\mu \geq \mu^{\star}$. 
Moreover, any local minima of $\boldsymbol{PMP}_{\mu}$ with $\mu > \mu^{\star}$ are also local minima of $\boldsymbol{MP}$.
\end{proposition}

%

Partial calmness has been studied for a range of problems, particularly bi-level problems \cite{ye1995,ye1997exact,mehlitz2020note}. Of specific importance to us is a related result due to Clarke \cite[Proposition 2.4.3]{clarke1990optimization} which allows for exact penalization (or equivalently partial calmness) and reads in our notation as follows.

\begin{proposition}\label{Proposition: LipschitzPenalty}
Consider the mathematical program $\boldsymbol{MP}$ and its penalized version $\boldsymbol{PMP}_{\mu}$. Assume that $f$ is Lipschitz continuous with Lipschitz constant $L$, the equality constraint takes the form of a distance to a closed set $S \subset C$,
\begin{equation*}
0 = h(x) = \text{distance}(x,S) = \text{inf}_{y \in S} \|x-y\| \,,
\end{equation*}
and $\boldsymbol{MP}$ attains a minimum. Then for $\mu > \mu^{\star} = L$, any local minimum of $\boldsymbol{PMP}_{\mu}$ is also a local minimum of $\boldsymbol{MP}$.
\end{proposition}

We note that $\|\cdot\|$ in Proposition~\ref{Proposition: LipschitzPenalty} can be an arbitrary norm. For a more general problem setup with a parametric set $S$ depicting the value function of an (inner) optimization problem, the reader is referred to \cite[Theorem 2.5]{ye1997exact}. Additionally, the setup can be extended to (squared) {\em merit functions} as penalty functions \cite[Theorems 2.6 and 2.9]{ye1997exact}. These generalize the notion of distance but are easier to formulate and compute.%
%


\bibliographystyle{IEEEtran}
\bibliography{refs}


\begin{IEEEbiography}[{\includegraphics[width=1in,height=1.25in,clip,keepaspectratio]{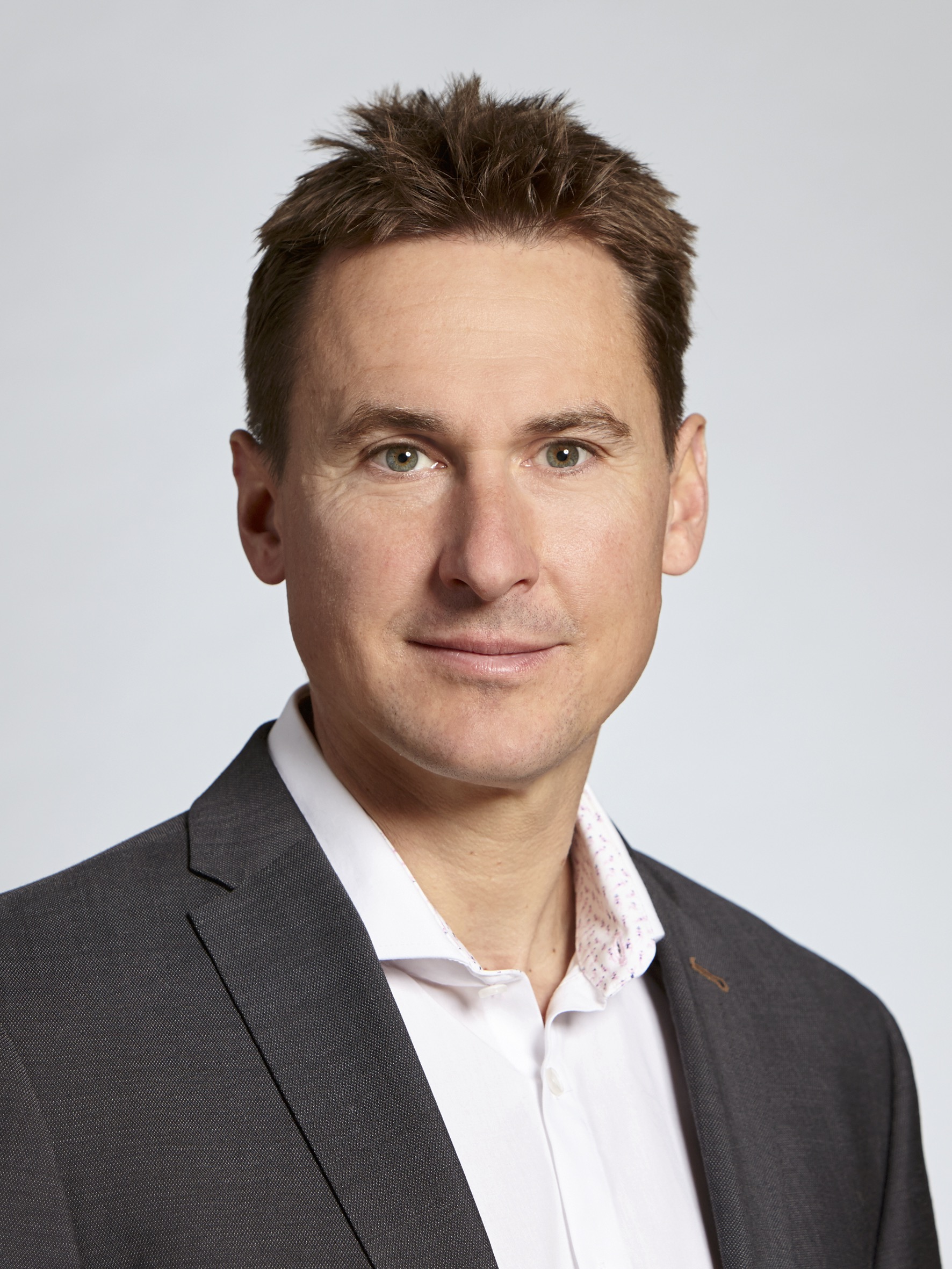}}]{Florian D\"orfler} (S'09-M'13-S'21) is an Associate Professor at the Automatic Control Laboratory at ETH Z\"urich and the Associate Head of the Department of Information Technology and Electrical Engineering. He received his Ph.D. degree in Mechanical Engineering from the University of California at Santa Barbara in 2013, and a Diplom degree in Engineering Cybernetics from the University of Stuttgart in 2008. From 2013 to 2014 he was an Assistant Professor at the University of California Los Angeles. His primary research interests are centered around control, optimization, and system theory with applications in network systems, in particular electric power grids. He is a recipient of the distinguished young research awards by IFAC (Manfred Thoma Medal 2020) and EUCA (European Control Award 2020). His students were winners or finalists for Best Student Paper awards at the European Control Conference (2013, 2019), the American Control Conference (2016), the Conference on Decision and Control (2020), the PES General Meeting (2020), and the PES PowerTech Conference (2017). He is furthermore a recipient of the 2010 ACC Student Best Paper Award, the 2011 O. Hugo Schuck Best Paper Award, the 2012-2014 Automatica Best Paper Award, the 2016 IEEE Circuits and Systems Guillemin-Cauer Best Paper Award, and the 2015 UCSB ME Best PhD award.
\end{IEEEbiography}

\begin{IEEEbiography}[{\includegraphics[width=1in,height=1.25in,clip,keepaspectratio]{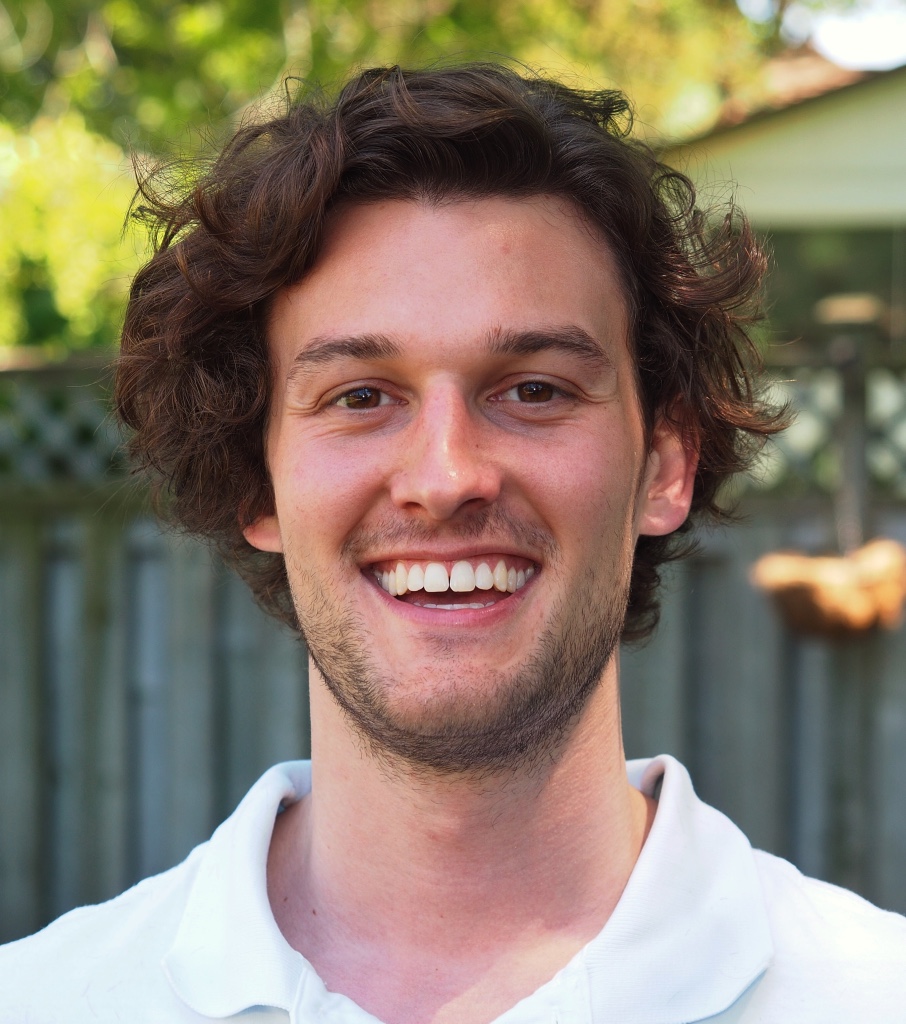}}]{Jeremy Coulson} (S'09-M'13)  is a PhD student with the Automatic Control Laboratory at ETH Z\"urich. He received his Master of Applied Science in Mathematics \& Engineering from Queen's University, Canada in August 2017. He received his B.Sc.Eng degree in Mechanical Engineering \& Applied Mathematics from Queen's University in 2015. His research interests include data-driven control methods, and stochastic optimization.
\end{IEEEbiography}

\begin{IEEEbiography}
[{\includegraphics[width=1in,height=1.25in,clip,keepaspectratio]{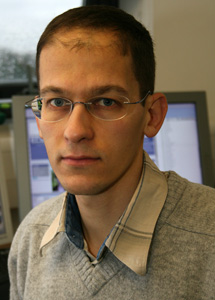}}]{Ivan Markovsky}
is an Associate Professor at the department ELEC of the Vrije Universiteit Brussel. He received his Ph.D. degree in Electrical Engineering from the Katholieke Universiteit Leuven in February 2005. From 2006 to 2012 he was an Assistant Professor at the School of Electronics and Computer Science of the University of Southampton. He is a recipient of an ERC starting grant "Structured low-rank approximation: Theory, algorithms, and applications" 2010--2015, Householder Prize honorable mention 2008, and research mandate by the Vrije Universiteit Brussel research council 2012--2022. His main research interests are computational methods for system theory, identification, and data-driven control in the behavioral setting.
\end{IEEEbiography}

\end{document}